\documentclass[11pt]{article}

\usepackage{amssymb}

\newcommand{\newc}{\newcommand}

\newc{\eqnoset}{\setcounter{equation}{0}}

\newcommand{\mref}[1]{(\ref{#1})}
\newcommand{\reflemm}[1]{Lemma~\ref{#1}}
\newcommand{\refrem}[1]{Remark~\ref{#1}}
\newcommand{\reftheo}[1]{Theorem~\ref{#1}}

\newcommand{\refcoro}[1]{Corollary~\ref{#1}}

\newcommand{\refsec}[1]{Section~\ref{#1}}

\newcommand{\beq}{\begin{equation}}
\newcommand{\eeq}{\end{equation}}
\newcommand{\beqno}[1]{\begin{equation}\label{#1}}

\newcommand{\barr}{\begin{array}}
\newcommand{\earr}{\end{array}}

\newc{\bearr}{\begin{eqnarray*}}
\newc{\eearr}{\end{eqnarray*}}

\newc{\bearrno}[1]{\begin{eqnarray}\label{#1}}
\newc{\eearrno}{\end{eqnarray}}

\newc{\non}{\nonumber}
\newc{\nol}{\nonumber\nl}

\newcommand{\bdes}{\begin{description}}
\newcommand{\edes}{\end{description}}
\newc{\benu}{\begin{enumerate}}
\newc{\eenu}{\end{enumerate}}
\newc{\btab}{\begin{tabular}}
\newc{\etab}{\end{tabular}}

\newtheorem{theorem}{Theorem}[section]
\newtheorem{defi}[theorem]{Definition}
\newtheorem{lemma}[theorem]{Lemma}
\newtheorem{rem}[theorem]{Remark}
\newtheorem{exam}[theorem]{Example}
\newtheorem{propo}[theorem]{Proposition}
\newtheorem{corol}[theorem]{Corollary}

\newcommand{\btheo}[1]{\begin{theorem}\label{#1}}
\newc{\brem}[1]{\begin{rem}\label{#1}\em}
\newc{\bexam}[1]{\begin{exam}\label{#1}\em}
\newc{\bdefi}[1]{\begin{defi}\label{#1}}
\newcommand{\blemm}[1]{\begin{lemma}\label{#1}}
\newcommand{\bprop}[1]{\begin{propo}\label{#1}}
\newcommand{\bcoro}[1]{\begin{corol}\label{#1}}
\newcommand{\etheo}{\end{theorem}}
\newcommand{\elemm}{\end{lemma}}
\newcommand{\eprop}{\end{propo}}
\newcommand{\ecoro}{\end{corol}}
\newc{\erem}{\end{rem}}
\newc{\eexam}{\end{exam}}
\newc{\edefi}{\end{defi}}

\newc{\rmk}[1]{{\bf REMARK #1: }}
\newc{\DN}[1]{{\bf DEFINITION #1: }}

\newcommand{\bproof}{{\bf Proof:~~}}
\newc{\eproof}{{\vrule height8pt width5pt depth0pt}\vspace{3mm}}

\newc{\bfrac}[2]{\dspl{\frac{#1}{#2}}}

\newc{\nid}{\noindent}

\newcommand{\dspl}{\displaystyle}
\newc{\grad}{\nabla}
\newc{\Div}{\mbox{div}}
\newc{\pdt}[1]{\dspl{\frac{\partial{#1}}{\partial t}}}
\newc{\pdn}[1]{\dspl{\frac{\partial{#1}}{\partial \nu}}}
\newc{\pdNi}[1]{\dspl{\frac{\partial{#1}}{\partial \mathcal{N}_i}}}
\newc{\pD}[2]{\dspl{\frac{\partial{#1}}{\partial #2}}}
\newc{\dt}{\dspl{\frac{d}{dt}}}
\newc{\bdry}[1]{\mbox{$\partial #1$}}
\newc{\sgn}{\mbox{sign}}

\newc{\Hess}[1]{\frac{\partial^2 #1}{\pdh z_i \pdh z_j}}
\newc{\hess}[1]{\partial^2 #1/\pdh z_i \pdh z_j}

\newc{\ag}{\alpha}
\newc{\bg}{\beta}
\newc{\cg}{\gamma}\newc{\Cg}{\Gamma}
\newc{\dg}{\delta}\newc{\Dg}{\Delta}
\newc{\eg}{\varepsilon}
\newc{\zg}{\zeta}
\newc{\thg}{\theta}
\newc{\llg}{\lambda}\newc{\LLg}{\Lambda}
\newc{\kg}{\kappa}
\newc{\rg}{\rho}
\newc{\sg}{\sigma}\newc{\Sg}{\Sigma}
\newc{\tg}{\tau}
\newc{\fg}{\phi}\newc{\Fg}{\Phi}
\newc{\vfg}{\varphi}
\newc{\og}{\omega}\newc{\Og}{\Omega}
\newc{\pdh}{\partial}

\newc{\ccG}{{\cal G}}

\newc{\ii}[1]{\int_{#1}}
\newc{\iidx}[2]{{\dspl\int_{#1}~#2~dx}}
\newc{\bii}[1]{{\dspl \ii{#1} }}
\newc{\biii}[2]{{\dspl \iii{#1}{#2} }}
\newc{\su}[2]{\sum_{#1}^{#2}}
\newc{\bsu}[2]{{\dspl \su{#1}{#2} }}

\newc{\biiom}[1]{{\dspl\int_{\bdrom}~ #1 ~d\sg}}
\newc{\io}[1]{{\dspl\int_{\Og}~ #1 ~dx}}
\newc{\bio}[1]{{\dspl\int_{\bdrom}~ #1 ~d\sg}}
\newc{\bsir}{\bsu{i=1}{r}}
\newc{\bsim}{\bsu{i=1}{m}}

\newc{\iibr}[2]{\iidx{\bprw{#1}}{#2}}
\newc{\Intbr}[1]{\iibr{R}{#1}}
\newc{\intbr}[1]{\iibr{\rg}{#1}}
\newc{\intt}[3]{\int_{#1}^{#2}\int_\Og~#3~dxdt}

\newc{\itQ}[2]{\dspl{\int\hspace{-2.5mm}\int_{#1}~#2~dz}}
\newc{\mitQ}[2]{\dspl{\rule[1mm]{4mm}{.3mm}\hspace{-5.3mm}\int\hspace{-2.5mm}\int_{#1}~#2~dz}}
\newc{\mitQQ}[3]{\dspl{\rule[1mm]{4mm}{.3mm}\hspace{-5.3mm}\int\hspace{-2.5mm}\int_{#1}~#2~#3}}

\newc{\mitx}[2]{\dspl{\rule[1mm]{3mm}{.3mm}\hspace{-4mm}\int_{#1}~#2~dx}}
\newc{\mitmu}[2]{\dspl{\rule[1mm]{3mm}{.3mm}\hspace{-4mm}\int_{#1}~#2~d\mu}}
\newc{\iidmu}[2]{{\dspl\int_{#1}~#2~d\mu}}

\newc{\iidm}[3]{{\dspl\int_{#1}~#2~d #3}}

\newc{\itQmu}[2]{\dspl{\int\hspace{-2.5mm}\int_{#1}~#2~d\mu}}
\newc{\mitQmu}[2]{\dspl{\rule[1mm]{4mm}{.3mm}\hspace{-5.3mm}\int\hspace{-2.5mm}\int_{#1}~#2~d\mu}}

\newc{\itQmuz}[2]{\dspl{\int\hspace{-2.5mm}\int_{#1}~#2~d\mu d\tau}}
\newc{\mitQmuz}[2]{\dspl{\rule[1mm]{4mm}{.3mm}\hspace{-5.3mm}\int\hspace{-2.5mm}\int_{#1}~#2~d\mu d\tau}}

\newc{\mitQq}[2]{\dspl{\rule[1mm]{4mm}{.3mm}\hspace{-5.3mm}\int\hspace{-2.5mm}\int_{#1}~#2~d\bar{z}}}
\newc{\itQq}[2]{\dspl{\int\hspace{-2.5mm}\int_{#1}~#2~d\bar{z}}}

\newc{\pder}[2]{\dspl{\frac{\partial #1}{\partial #2}}}

\newc{\bdrom}{\bdry{\Og}}

\newc{\bilhom}{\mbox{Bil}(\mbox{Hom}(\RR^{nm},\RR^{nm}))}
\newc{\VV}[1]{{V(Q_{#1})}}

\newc{\ccA}{{\mathcal A}}
\newc{\ccB}{{\mathcal B}}
\newc{\ccC}{{\mathcal C}}
\newc{\ccD}{{\mathcal D}}
\newc{\ccE}{{\mathcal E}}
\newc{\ccH}{\mathcal{H}}
\newc{\ccF}{\mathcal{F}}
\newc{\ccI}{{\mathcal I}}
\newc{\ccJ}{{\mathcal J}}
\newc{\ccK}{{\mathcal K}}
\newc{\ccP}{{\mathcal P}}
\newc{\ccQ}{{\mathcal Q}}
\newc{\ccR}{{\mathcal R}}
\newc{\ccS}{{\mathcal S}}
\newc{\ccT}{{\mathcal T}}
\newc{\ccX}{{\mathcal X}}
\newc{\ccY}{{\mathcal Y}}
\newc{\ccZ}{{\mathcal Z}}

\newc{\bb}[1]{{\mathbf #1}}

\newc{\myprod}[1]{\langle #1 \rangle}
\newc{\mypar}[1]{\left( #1 \right)}

\newc{\BLLg}{\mathbf{\LLg}}

\newc{\mA}{\mathbf{A}}
\newc{\mB}{\mathbf{B}}
\newc{\mC}{\mathbf{C}}
\newc{\mD}{\mathbf{D}}
\newc{\mE}{\mathbf{E}}
\newc{\mF}{\mathbf{F}}
\newc{\mJ}{\mathbf{J}}
\newc{\mG}{\mathbf{G}}
\newc{\mP}{\mathbf{P}}
\newc{\mR}{\mathbf{R}}
\newc{\mQ}{\mathbf{Q}}
\newc{\mX}{\mathbf{X}}
\newc{\muu}{\mathbf{u}}
\newc{\mvv}{\mathbf{v}}

\newc{\mllg}{\mathbb{\lambda}}
\newc{\mLLg}{\mathbf{\LLg}}

\newc{\lspn}[2]{\mbox{$\| #1\|_{\Lsp{#2}}$}}
\newc{\Lpn}[2]{\mbox{$\| #1\|_{#2}$}}
\newc{\Hn}[1]{\mbox{$\| #1\|_{H^1(\Og)}$}}

\newc{\mynorm}[2]{\| #1\|_{#2}}

\newcommand{\RR}{{\rm I\kern -1.6pt{\rm R}}}

\newc{\itQQ}[2]{\dspl{\int_{#1}#2\,dz}}
\newc{\mmitQQ}[2]{\dspl{\rule[1mm]{4mm}{.3mm}\hspace{-4.3mm}\int_{#1}~#2~dz}}
\newc{\MmitQQ}[2]{\dspl{\rule[1mm]{4mm}{.3mm}\hspace{-4.3mm}\int_{#1}~#2~d\mu}}

\newc{\MUmitQQ}[3]{\dspl{\rule[1mm]{4mm}{.3mm}\hspace{-4.3mm}\int_{#1}~#2~d#3}}
\newc{\MUitQQ}[3]{\dspl{\int_{#1}~#2~d#3}}

\newc{\mccP}{\mathbb{P}}
\newc{\mccK}{\mathbb{K}}

\newc{\DKTmU}{\mccK(U)}
\newc{\DKTmUold}{(K_U(U)^{-1})^T}

\newc{\myPi}{\mathbf{W}}
\newc{\myIbar}{\bar{\ccI}_1}
\newc{\myIhat}{\hat{\ccI}_1}
\newc{\myIbreve}{\breve{\ccI}_0}

\newc{\mmk}{\mathbf{k}}

\newc{\extraI}{\mathbb{I}}
\newc{\mccIe}{\ccI}
\newc{\mccIee}{\ccF}

\newc{\mU}{\mathbf{U}}
\newc{\mL}{\mathbf{L}}

\begin{document}

\vspace*{-.8in}
\begin{center} {\LARGE\em On the Solvability of a Class of Degenerate or Singular Strongly Coupled Parabolic Systems.}

 \end{center}

\vspace{.1in}

\begin{center}

{\sc Dung Le}{\footnote {Department of Mathematics, University of
Texas at San
Antonio, One UTSA Circle, San Antonio, TX 78249. {\tt Email: Dung.Le@utsa.edu}
}}

\end{center}

\begin{abstract}
The existence of strong solutions to  general class of strongly coupled parabolic systems will be discussed. These systems can be degenerate or singular as boundedness of theirs solutions are unavailable and not assummed. The results greatly improve those in a recent papers \cite{letrans,dleJFA,dleANS} as the systems can have quadratic growth in gradients. A unified proof for both cases is presented. Most importantly, the VMO assumption in \cite{dleJFA,dleANS} will be replaced by a much versatile one thanks to a new local weighted Gagliardo-Nirenberg involving BMO norms. Degenerate and singular generalized SKT models in biology will be presented as a nontrivial application of the main theorem. \end{abstract}

\vspace{.2in}

\section{Introduction}\label{introsec}\eqnoset

In this paper, for any $T_0>0$ and bounded domain $\Og$ with smooth boundary in $\RR^n$, $n\ge2$, we consider the following parabolic  system of $m$ equations ($m\ge2$) for the unknown $u:Q\to \RR^m$, where $Q=\Og\times(0,T_0)$
\beqno{e1}\left\{\barr{l} u_t-\Div(A(x,u)Du)=\hat{f}(x,u,Du),\quad (x,t)\in Q,,\\\mbox{$u=0$ or $\frac{\partial u}{\partial \nu}=0$ on $\partial \Og$}, \, t\in (0,T_0),\\u(x,0)=U_0(x) \quad x\in \Og. \earr\right.\eeq
Here, $A(x,u)$ is a $m\times m$ matrix in $x\in\Og$ and $u\in\RR^m$, $\hat{f}:\Og\times\RR^m\times \RR^{mn}\to\RR^m$ is a vector valued function. 
The initial data $U_0$ is given in $W^{1,p_0}(\Og,\RR^m)$ for some $p_0>n$, the dimension of $\Og$. As usual, $W^{1,p}(\Og,\RR^m)$, $p\ge1$, will denote the standard Sobolev spaces whose elements are vector valued functions $u\,:\,\Og\to \RR^m$ with finite norm $$\|u\|_{W^{1,p}(\Og,\RR^m)} = \|u\|_{L^p(\Og)} + \|Du\|_{L^p(\Og)}.$$

The strongly coupled system \mref{e1} appears in many physical applications, for instance, Maxwell-Stephan systems describing the diffusive transport of multicomponent mixtures, models in reaction and diffusion in electrolysis, flows in porous media, diffusion of polymers, or population dynamics, among others.

We will discuss the existence of strong solutions to \mref{e1}. We say that $u$ is a strong solution if $u$ solves \mref{e1} a.e. on $\bar{Q}$ with $Du\in L^\infty_{loc}(Q)$ and $D^2u\in L^2_{loc}(Q)$.

It is always assumed that the matrix $A(x,u)$ is elliptic in the sense that there exist two scalar positive continuous functions $\llg_1(x,u), \llg_2(x,u)$ such that \beqno{genelliptic} \llg_1(x,u)|\zeta|^2 \le \myprod{A(x,u)\zeta,\zeta}\le \llg_2(x,u)|\zeta|^2 \quad \mbox{for all } x\in\Og,\,u\in\RR^m,\,\zeta \in\RR^{nm}.\eeq If there exist positive constants $c_1,c_2$ such that $c_1\le \llg_1(x,u)$ and  $\llg_2(x,u)\le c_2$ then we say that $A(x,u)$ is {\em regular elliptic}. If $c_1\le \llg_1(x,u)$ and  $\llg_2(x,u)/\llg_1(x,u)\le c_2$, we say that $A(x,u)$ is {\em uniform elliptic}. On the other hand, if we allow $c_1=0$ and $\llg_1(x,u)$ tend to zero (respectively, $\infty$) when $|u|\to\infty$ then we say that $A(x,u)$ is {\em singular} (respectively, {\em degenerate}). 

We consider the following structural conditions on the data of \mref{e1}.

\bdes

\item[A)] $A(x,u)$ is $C^1$ in $x\in\Og$, $u\in\RR^m$ and there exist a constant $C_*>0$  and  scalar $C^1$ positive functions $\llg(u),\og(x)$  such that  for all $u\in\RR^m$,  $\zeta\in\RR^{mn}$ and $x\in\Og$ 
\beqno{A1} \llg(u)\og(x)|\zeta|^2 \le \myprod{A(x,u)\zeta,\zeta} \mbox{ and } |A(x,u)|\le C_*\llg(u)\og(x).\eeq 

In addition, there is a constant $C$ such that $|\llg_u(u)||u|\le C\llg(u)$ and  \beqno{Axcond} |A_u(x,u)|\le C|\llg_u(u)|\og(x),\; |A_x(x,u)|\le C|\llg(u)||D\og|.\eeq

\edes

Here and throughout this paper, if $B$ is a $C^1$ (vector valued) function in $u\in \RR^m$ then we abbreviate its  derivative $\frac{\partial B}{\partial u}$ by $B_u$. Also, with a slight abuse of notations, $A(x,u)\zeta$, $\myprod{A(x,u)\zeta,\zeta}$ in \mref{genelliptic}, \mref{A1} should be understood in the following way: For $A(x,u)=[a_{ij}(x,u)]$, $\zeta\in\RR^{mn}$ we write $\zeta=[\zeta_i]_{i=1}^m$ with $\zeta_i=(\zeta_{i,1},\ldots\zeta_{i,n})$ and $$A(x,u)\zeta=[\Sigma_{j=1}^m a_{ij}\zeta_j]_{i=1}^m,\; \myprod{A(x,u)\zeta,\zeta}=\Sigma_{i,j=1}^m a_{ij}\myprod{\zeta_i,\zeta_j}.$$

We also assume that $A(x,u)$ is regular elliptic for {\em bounded} $u$.

\bdes\item[AR)] $\og\in C^1(\Og)$ and there are positive numbers $\mu_*,\mu_{**}$ such that 
\beqno{Aregcond2} \mu_*\le\og(x)\le \mu_{**}, \; |D\og(x)|\le \mu_{**} \quad \forall x\in\Og.\eeq For any bounded set $K\subset\RR^m$ there is a constant $\llg_*(K)>0$ such that
\beqno{Aregcond1} \llg_*(K)\le\llg(u) \quad \forall u\in K.\eeq \edes

Concerning the reaction term $\hat{f}(x,u,Du)$, which may have linear or {\em quadratic} growth in $Du$, we assume the following condition.
\bdes \item[F)] There exist a constant $C$ and a nonegative differentiable function $f:\RR^m\to\RR$ such that $\hat{f}$ satisfies: 
\beqno{specfcond}f(u)\le C|f_u(u)|(1+|u|).\eeq

For  any diffrentiable vector valued functions $u:\RR^n\to\RR^m$ and $p:\RR^n\to\RR^{mn}$ we assume either that
\bdes\item[f.1)] $\hat{f}$ has a linear growth in $p$ 
\beqno{FUDU11}|\hat{f}(x,u,p)| \le C\llg(u)|p|\og(x) + f(u)\og(x),\eeq $$|D\hat{f}(x,u,p)| \le  C(\llg(u)|Dp|+ |\llg_u(u)||p|^2)\og+ C\llg(u)|p||D\og|+C|D(f(u)\og(x))|;$$\edes or \bdes 
\item[f.2)] $\llg_{uu}(u)$ exists and $\hat{f}$ has a quadratic growth in $p$ \beqno{FUDU112}|\hat{f}(x,u,p)| \le C|\llg_u(u)||p|^2\og(x) + f(u)\og(x),\eeq  $$\barr{lll}|D\hat{f}(x,u,p)| &\le&  C(|\llg_u(u)||p||Dp|+ |\llg_{uu}(u)||p|^3)\og+ C|\llg_u(u)||p|^2|D\og|\\&&+C|D(f(u)\og(x))|.\earr$$ Furthermore, we assume that \beqno{llgquadcond} |\llg_{uu}(u)|\llg(u)\le C|\llg_u(u)|^2.\eeq
\edes
\edes

By a formal differentiation of \mref{FUDU11} and \mref{FUDU112}, one can see that  the growth conditions for $\hat{f}$ naturally implies those of  $D\hat{f}$ in the above assumption. The condition \mref{llgquadcond} is verified easily if $\llg(u)$ has a polynomial growth in $|u|$.

The first fundamental problem in the study of \mref{e1} is the local and global existence of its solutions. One can decide to work with either weak or strong solutions. In the first case, the existence of a weak solution can be achieved via Galerkin, time discretization  or variational methods but its regularity (e.g., boundedness, H\"older continuity of the solution and its higher derivatives) is still an open issue. Several works have been done along this line to improve the early work \cite{GiaS} of Giaquinta and Struwe and establish {\em partial regularity} of {\em bounded} weak solutions to \mref{e1}.

Otherwise, if strong solutions are considered then theirs existence can be established via semigroup theories as in the works of Amann \cite{Am1,Am2}. Combining with interpolation theories of Sobolev's spaces, Amann established local and global existence of a strong solution $u$ of \mref{e1} under the assumption that one can controll $\|u\|_{W^{1,p}(\Og,\RR^m)}$ for some $p>n$. His theory did not apply to the case where $\hat{f}$ has quadratic growth in $Du$ as in f.2).

In both forementioned approaches, the assumption on the boundedness of $u$ must be the starting point and the techniques in both cases  rely heavily on the fact that $A(x,u)$ is {\em regular elliptic}. 
For strongly coupled systems like \mref{e1}, as invariant/maximum principles for cross diffusion systems are generally unavailable, the boundedness of the solutions is already a hard problem. One usually needs to use ad hoc techniques on the case by case basis to show that $u$ is bounded (see \cite{kuf1,Red}).  Even for bounded weak solutions, we know that they are only H\"older continuous almost everywhere (see \cite{GiaS}).
In addition, there are counter examples for systems ($m>1$) which exhibit solutions that start smoothly and remain bounded but develop singularities in higher norms in finite times (see \cite{JS}).

In our recent work \cite{dleANS,dleJFA}, we choose a different approach making use of fixed point theory and discussing the existence of {\em strong} solutions of \mref{e1} under the weakest assumption that they are a-priori VMO, not necessarily bounded, and general structural conditions on the data of \mref{e1} which are independent of $x$, we assumed only that $A(u)$ is {\em uniformly elliptic}.  Applications were presented in \cite{dleJFA} when $\llg(u)$ has a positive polynomial growth in $|u|$ and, without the boundedness assumption on the solutions, so \mref{e1} can be degenerate as $|u|\to\infty$. The singular case, $\llg(u)\to 0$ as $|u|\to\infty$, was not discussed there.

In this paper, we will establish much stronger results than those in \cite{dleANS} under much more general assumptions on the structure of \mref{e1} as described in A) and F). Beside the minor fact that the data can depend on $x$, we allow further that:
\begin{itemize}
	\item $A(x,u)$ can be either degenerate or singular as $|u|$ tends to infinity;
	\item  $\hat{f}(x,u,Du)$ can have a quadratic growth in $Du$ as in f.2); \item no a-priori boundedness of solutions is assummed  but a a very weak integrability of strong solutions of \mref{e1} is considered.
	
\end{itemize}

Most remarkably, the key assumption in \cite{dleJFA,dleANS} that the BMO norm of $u$ is small in small balls will be replaced by a more versatile one in this paper: $K(u)$ is has small BMO norm in small balls for some suitable map $K:\RR^m\to\RR^m$. This allows us to consider the singular case where one may not be able to estimate the BMO norm of $u$  but that of $K(u)$. Examples of this case in applications will be provided in \refsec{res} where $|K(u)|\sim \log(|u|)$.

One of the key ingredients in the proof in \cite{dleANS,dleJFA} is the {\em local} weighted Gagliardo-Nirenberg inequality involving BMO norm \cite[Lemma 2.4]{dleANS}.  In this paper, we make use of a new version of this inequality reported in our work \cite{dleGNnew} replacing the BMO norm of $u$ by that of $K(u)$ for some suitable map $K:\RR^m\to\RR^m$.

We organize our paper as follows. In \refsec{res} we state the main result, \reftheo{gentheo1}, of this paper and its application to the generalized SKT systems on planar domain. In \refsec{GNsec} we recall the new version of the local weighted Gagliardo-Nirenberg inequality in \cite{dleGNnew} to prepare for the proof the main \reftheo{gentheo1} in \refsec{w12est}. The proof of solvability of the generalized SKT systems  in \refsec{res} is provided in \refsec{genSKTsec}.

\section{Preliminaries and Main Results}\eqnoset\label{res}

We state the main results of this paper in this section. The key assumption of these results is some uniform a priori estimate for the BMO norm of $K(u)$ where $K$ is some suitable map on $\RR^m$ and $u$ is  any  strong solution to \mref{e1}. To begin, we recall some basic definition in Harmonic Analysis.

Let $\og\in L^1(\Og)$ be a nonnegative function and  define the measure $d\mu=\og(x)dx$. For any $\mu$-measurable subset $A$ of $\Og$  and any  locally $\mu$-integrable function $U:\Og\to\RR^m$ we denote by  $\mu(A)$ the measure of $A$ and $U_A$ the average of $U$ over $A$. That is, $$U_A=\mitmu{A}{U(x)} =\frac{1}{\mu(A)}\iidmu{A}{U(x)}.$$
We define the measure $d\mu =\og(x)dx$ and recall that a vector valued function $f\in L^1(\Og,\mu)$ is said to be in $BMO(\Og,\mu)$ if \beqno{bmodef} [f]_{*,\mu}:=\sup_{B_R\subset\Og}\mitmu{B_R}{|f-f_{B_R}|}<\infty,\quad f_{B_R}:=\frac{1}{\mu(B_R)}\iidmu{B_R}{f}.\eeq We then define $$\|f\|_{BMO(\Og,\mu)}:=[f]_{*,\mu}+\|f\|_{L^1(\Og,\mu)}.$$

For $\cg\in(1,\infty)$ we say that a nonnegative locally integrable function $w$ belongs to the class $A_\cg$ or $w$ is an $A_\cg$ weight on $\Og$ if the quantity
\beqno{aweight} [w]_{\cg,\Og} := \sup_{B\subset\Og} \left(\mitmu{B}{w}\right) \left(\mitmu{B}{w^{1-\cg'}}\right)^{\cg-1} \quad\mbox{is finite}.\eeq
Here, $\cg'=\cg/(\cg-1)$. For more details on these classes we refer the reader to \cite{OP,st}. If the domain $\Og$ is specified we simply denote $ [w]_{\cg,\Og}$ by $ [w]_{\cg}$.

Throughout this paper, in  our statements and proofs, we use $C,C_1,\ldots$ to denote various constants which can change from line to line but depend only on the parameters of the hypotheses in an obvious way. We will write $C(a,b,\ldots)$ when the dependence of  a constant $C$ on its parameters is needed to emphasize that $C$ is bounded in terms of its parameters. We also write $a\lesssim b$ if there is a universal constant $C$ such that $a\le Cb$. In the same way, $a\sim b$ means $a\lesssim b$ and $b\lesssim a$.

To begin, as in \cite{dleANS} with $A$ is independent of $x$, we assume that the eigenvalues of the matrix $A(x,u)$ are not too far apart. Namely, for $C_*$ defined in \mref{A1} of A) we assume
\bdes\item[SG)] $(n-2)/n <C_*^{-1}$. \edes Here $C_*$ is, in certain sense, the ratio of the largest and smallest eigenvalues of $A(x,u)$. This condition seems to be necessary as we deal with systems, cf.  \cite{dleN}. 

First of all, we will assume that the system \mref{gensys} satisfies the structural conditions A) and F). Additional assumptions serving the purpose of this paper then follow  so that the local weighted Gagliardo-Nirenberg inequality of \cite{dleGNnew} can applies here.

\bdes 
\item[H)] There is a $C^1$ map $K:\RR^m\to\RR^m$ such that $\mathbb{K}(u)=(K_u(u)^{-1})^T$ exists and $\mathbb{K}_u\in L^\infty(\RR^m)$.
Furthermore,  for all $u\in \RR^m$
\beqno{kappamainz} |\mathbb{K}(u)|\lesssim \llg(u)|\llg_u(u)|^{-1}.\eeq
\edes

We consider the following system
\beqno{gensysmain}\left\{\barr{l}u_t -\Div(A(x,u)Du)=\hat{f}(x,u,Du),\quad x\in \Og,\\\mbox{$u=0$ or $\frac{\partial u}{\partial \nu}=0$ on $\partial \Og$}. \earr\right.\eeq

We imbed this system in the following family of systems 
\beqno{gensysfammain}\left\{\barr{l}u_t -\Div(A(x,\sg u)Du)=\hat{f}(x,\sg u,\sg Du),\quad x\in \Og, \sg\in[0,1],\\\mbox{$u=0$ or $\frac{\partial u}{\partial \nu}=0$ on $\partial \Og$}. \earr\right.\eeq

For any strong solution $u$ of \mref{gensysfammain} we will consider  the following assumptions.

\bdes 

\item[M.0)] There exists a constant $C_0$ such that for some $r_0>1$ and $\bg_0\in(0,1)$ \beqno{llgmainhyp} \sup_{\tau\in(0,T_0)}\||f_u(\sg u)|\llg^{-1}(\sg u)\|_{L^{r_0}(\Og,\mu)},\; \sup_{\tau\in(0,T_0)}\|u^{\bg_0}\|_{L^{1}(\Og,\mu)}\le C_0 ,\eeq
\beqno{hypoiterpis1}\itQmuz{Q}{(|f_u(\sg u)|+\llg(\sg u))(|Du|^{2}+|u|^2)}\le C_0.\eeq

\item[M.1)] For any given  
$\mu_0>0$ 
there is positive $R_{\mu_0}$ sufficiently small in terms of the constants in A) and F) such that \beqno{Keymu0} \sup_{x_0\in\bar{\Og},\tau\in(0,T_0)}\|K(\sg u)\|_{BMO(B_{R}(x_0)\cap\Og,\mu)}^2 \le \mu_0.\eeq  

Furthermore, for  $\myPi_p(\sg,x,\tau):= \llg^{p+\frac12}(\sg u)|\llg_u(\sg u)|^{-p}$ and any $p\in[1,n/2]$ there exist some $\ag>2/(p+2)$, $\bg<p/(p+2)$  such that $\sup_{\tau\in(0,T_0)}[\myPi_p^{\ag}]_{\bg+1,B_{R_{\mu_0}}(x_0)\cap\Og}\le C_0$. 
\edes

The main theorem of this paper is the following. 

\btheo{gentheo1} Assume A), F), AR) and H).  Moreover, if $\hat{f}$ has a quadratic growth in $Du$ as in f.2) then we assume also that $n\le3$. Suppose also that any strong solution $u$ to \mref{gensysfammain} satisfies M.0), M.1) 
uniformly in $\sg\in[0,1]$.  

Then the system \mref{gensysmain} has a unique strong solution on $\Og\times(0,T_0)$. \etheo

The condition \mref{Keymu0} on the smallness of the BMO norm of $K(u)$ in small balls is the most crucial one in applications. In \cite{dleANS, dleJFA}, we consider the case $\llg(u)\sim(\llg_0+|u|)^k$ with $k>0$ and assume that $K=Id$, the identity matrix. We assumed that $K(u)=u$ has small BMO norm in small balls, which can be verified by establishing that $\|Du\|_{L^n(\Og)}$ is bounded. These results already improve those of Amann in \cite{Am1,Am2} where boundedness of solutions was assumed and uniform estimates for $\|Du\|_{L^p(\Og)}$ for some $p>n$ is needed. Both of such conditions seems to be very difficult to be verified in applications.

We should remark that all the assumptions on strong solutions of the family \mref{gensysfammain} can be checked by considering the case $\sg=1$ (i.e. \mref{e1}) because these systems satisfy the same structural conditions uniformly with respect to the parameter $\sg\in[0,1]$.

We present an application of \reftheo{gentheo1}. This example concerns cross diffusion systems with polynomial growth data on planar domains. This type of systems occurs in many applications in mathematical biology and ecology. An famous example of such systems is the SKT model (see \cite{dleJFA,SKT,yag}) for two species with population densities $u,v$ satifying \beqno{SKTor} \left\{\barr{lll}u_t&=&\Delta(u[d_1+\ag_{11}u+\ag_{12}v])+f_1(u,v),\\v_t&=&\Delta(v[d_2+\ag_{21}u+\ag_{22}v])+f_2(u,v).\earr\right.\eeq

We consider the following generalized SKT system  with Dirichlet or Neumann boundary conditions on a bounded domain $\Og\subset\RR^n$ for vector valued unknown $u:\Og\times(0,T_0)\to \RR^m$.
\beqno{genSKT} u_t-\Delta(P_i(u))=B_i(u,Du)+ f_i(u), \quad i=1,\ldots,m.\eeq Here, $P_i:\RR^m\to\RR$ are $C^2$ functions. The functions $B_i, f_i$ are $C^1$ functions on $\RR^m\times\RR^{mn}$ and $\RR^m$ respectively. We will assume that $B_i(u,Du)$ has the following linear growth in $Du$. 
\beqno{Bgrowth} |B_i(u,Du)|\lesssim \llg^\frac12(u)|Du|.\eeq  The growth in $Du$ of $B_i(u,Du)$ is a bit different from f.1) in this paper but we will see that \reftheo{gentheo1} still applies here (see \refrem{SKTgenBrem}).

The system \mref{genSKT} generalizes \mref{SKTor} by letting $P_i(u)=u_i\llg_i(u)$ for some $C^2$ functions $\llg_i(u)$ and consider the following assumption (see also \refrem{SKTgenrem} after the theorem). 

\bdes\item[L)] There exist $C^2$ nonnegative scalar functions $\llg_i$, $i=1,\ldots,m$, and $\llg$ on $\RR^m$ such that \beqno{llgcondz} \llg_i(u)\lesssim\llg(u),\; |u||\llg_u(u)|\lesssim \llg(u),\;|(\llg_i(u))_{uu}|\lesssim |\llg_{uu}(u)|.\eeq

The matrices $\mL=\mbox{diag}[\llg_1(u),\ldots,\llg_m(u)]$ and, with a slight abuse of notation, $\mL_u=D_u[\llg_i(u)]_{i=1}^m$ satisfy the following conditions. 
\beqno{Aellipticmain} \myprod{(\mL+\mbox{diag}[u_1,\ldots,u_m]\mL_u)\zeta,\zeta}\ge \llg(u)|\zeta|^2,\eeq
\beqno{mLhypmain} |\mL_u|\lesssim |\llg_u(u)|,\; |\mL_u^{-1}|\lesssim |\llg_u(u)|^{-1}.\eeq \edes

As $\Delta(P_i(u))=\Div(A(u)Du)$ with $A(u)=(\mL+\mbox{diag}[u_i]\mL_u)$, the condition \mref{Aelliptic} is necessary for \mref{genSKT} being elliptic. In fact, if $|(\llg_i(u)_u||u|\le c_i\llg_i(u)$ for some small $c_i$ then it is not difficult to see that \mref{Aelliptic} holds.

We now embed \mref{genSKT} into the following family of system \beqno{genSKTfam} u_t-\Div(A(\sg u)Du)=B_i(\sg u,\sg Du)+ f_i(\sg u), \quad \sg\in[0,1],\,i=1,\ldots,m.\eeq

As a consequence of \reftheo{gentheo1}, we will have the following.
\btheo{SKTgenthm} Assume L) and \mref{Bgrowth}. Assume further that $n=2$, $\llg$ satisfies AR) and there is a constant $C_0$ such that \beqno{LLubound0} |\llg_u(u)|\llg^{-2}(u)\le C_0 \mbox{ for all $u\in\RR^m$}.\eeq In addition to the integrability condition M.0), assume that any strong solution $u$ to \mref{genSKTfam} satisfies  \beqno{fSKT2d0} \dspl{\int_0^{T_0}}\iidx{\Og}{|\llg(\sg u)f_i(\sg u)|^2}\le C_0 \mbox{ for all $i=1,\ldots,m$}.\eeq Then \mref{genSKT} has a unique strong solution on $\Og\times(0,T_0)$.

\etheo

If $\llg(u)\sim (\llg_0+|u|)^k$ then \mref{LLubound0} holds if $k>-1$. Therefore, this theorem inludes the singular case of SKT system when $|u|$ becomes unbounded. 

\brem{SKTgenrem} The condition \mref{mLhypmain} is inspired by the SKT system \mref{SKTor}. In fact,
let $\ag_i=[\ag_{ij}]_{j=1}^m$ be $m$ linearly independent vectors in $\RR^m$. For some $k>0$ and $d_i>0$ we define $\llg_i(u)=d_i+\myprod{u,\ag_i}^k$ and $\ag:=[\ag_i]_{i=1}^m$. Then $\partial_u\llg_i(u)= k\myprod{u,\ag_i}^{k-1}\ag_i^T$ so that $\mL_u=k\mbox{diag}[\myprod{u,\ag_i}^{k-1}]\ag$ and
$\mL_u^{-1}=k^{-1}\ag^{-1}\mbox{diag}[\myprod{u,\ag_i}^{-k+1}]$. If $\myprod{u,\ag_i}\sim\myprod{u,\ag_j}$ for $i\ne j$ then $|\mL_u^{-1}|\sim |\ag^{-1}||\llg_u(u)|^{-1}$ with $\llg(u)=\sum_i(d_i+|\myprod{u,\ag_i}|^k)$. The system \mref{genSKT} is degenerate when $|u|\to\infty$. We see that the SKT system \mref{SKTor} is included in this case for $m=2,k=1$. 

On the other hand, we can consider the singular case when $k<0$. We define $\llg_i(u)=(d_i+\myprod{u,\ag_i})^k$. Then $\partial_u\llg_i(u)= k(d_i+\myprod{u,\ag_i})^{k-1}\ag_i^T$ so that $\mL_u=k\mbox{diag}[(d_i+\myprod{u,\ag_i})^{k-1}]\ag$ and
$\mL_u^{-1}=k^{-1}\ag^{-1}\mbox{diag}[(d_i+\myprod{u,\ag_i})]$. We then have $|\mL_u^{-1}|\sim |\ag^{-1}||\llg_u(u)|^{-1}$ with $\llg(u)=\sum_i|d_i+\myprod{u,\ag_i}|^k$. In both cases, we see that \mref{mLhypmain} holds. \erem

\section{A general  local weighted Gagliardo-Nirenberg inequality} \eqnoset\label{GNsec}
In this section, we present a  local weighted Gagliardo-Nirenberg inequality in our recent work \cite{dleGNnew}, which will be one of the main ingredients of the proof of our main technical theorem in \refsec{w12est}. 
This inequality generalizes \cite[Lemma 2.4]{dleANS} by replacing the Lebesgue measure with general one and the BMO norm of $u$ with that of $K(u)$ where $K$ is a suitable map on $\RR^m$, and so the applications of our main technical theorem in the next section will be much more versatile than those in \cite{dleJFA,dleANS}.

Let us begin by describing the assumptions in \cite{dleGNnew} for this general inequality. We say that $\Og$ and $\mu$ support a $q_*$-Poincar\'e inequality if the following holds. \bdes \item[P)]  There exist $q_*\in(0,2]$, $\tau_*\ge 1$ and some constant $C_P$ such that \beqno{Pineq2} 
\mitmu{B}{|h-h_{B}|}\le C_Pl(B)
\left(\mitmu{\tau_*B}{|Dh|^{q_*}}\right)^\frac{1}{q_*}\eeq
for  any cube $B\subset\Og$ with side length $l(B)$ and any function $u\in C^1(B)$. \edes

Here and throughout this section,  we denote by $l(B)$ the side length of $B$ and by $\tau B$ the cube which is concentric with $B$ and has side length $\tau l(B)$. We also write $B_R(x)$ for a cube centered at $x$ with side length $R$ and sides parallel to to standard axes of $\RR^n$. We will omit $x$ in the notation $B_R(x)$ if no ambiguity can arise.

We consider the following conditions on the measure $\mu:=\og(x)dx$ for the validity of \mref{Pineq2} (see \cite{Haj}).

\bdes \item[LM.1)]  For some $N\in(0,n]$ and any ball $B_r$ we have $\mu(B_r)\le C_\mu r^N$. 
Assume also that $\mu$ supports the 2-Poincar\'e inequality \mref{Pineq2} in P). Furthermore, $\mu$ is doubling and satisfies  the following inequality for some $s_*>0$ \beqno{fracmu} \left(\frac{r}{r_0}\right)^{s_*}\le C_\mu\frac{\mu(B_r(x))}{\mu(B_{r_0}(x_0))},\eeq where $B_r(x), B_{r_0}(x_0)$ are any cubes with $x\in B_{r_0}(x_0)$.

\item[LM.2)]  $\og=\og_0^2$ for some $\og_0\in C^1(\Og)$ and $d\mu=\og_0^2 dx$  also supports a Hardy type inequality: There is a constant $C_H$ such that for any function $u\in C^1_0(B)$\beqno{lehr1m} \iidx{\Og}{|u|^2|D\og_0|^2}\le C_H\iidx{\Og}{|Du|^2\og_0^2}.\eeq 
\edes

We assume the following hypotheses.

\bdes
\item[A.1)] Let $K:\mbox{dom}(K)\to\RR^m$ be a $C^1$ map on a domain $\mbox{dom}(K)\subset\RR^m$ such that $\DKTmU=(K_U(U)^{-1})^T$ exists and $\mccK_U\in L^\infty(\mbox{dom}(K))$.

Furthermore, let $\Fg,\LLg:\mbox{dom}(K)\to\RR^+$ be $C^1$ positive functions. We assume that for all $U\in \mbox{dom}(K)$
\beqno{kappamain} |\DKTmU|\lesssim \LLg(U)\Fg^{-1}(U),\eeq
\beqno{logcondszmain}|\Fg_U(U)||\mccK(U)|\lesssim \Fg(U).\eeq
\edes

Let $\Og_*$ be a proper subset of $\Og$ and $\og_*$ be a function in $C^1(\Og)$ satisfying  
\beqno{subogm} \og_*\equiv 1 \mbox{ in $\Og_*$ and } \og_*\le 1 \mbox{ in $\Og$}.\eeq

For any $U\in C^2(\Og,\mbox{dom}(K))$  we denote
\beqno{Idefm} I_1:=\iidmu{\Og}{\Fg^2(U)|DU|^{2p+2}},\;
I_2:=\iidmu{\Og}{\LLg^2(U)|DU|^{2p-2}|D^2U|^2},\eeq
\beqno{Idef1zm} \myIbar:=\iidmu{\Og}{|\LLg_U(U)|^2|DU|^{2p+2}},\;I_{1,*}:=\iidmu{\Og_*}{\Fg^2(U)|DU|^{2p+2}},\eeq \beqno{I0*m} \breve{I}_{0,*}:=\sup_\Og|D\og_*|^2\iidmu{\Og}{\LLg^2(U)|DU|^{2p}}.\eeq 

We established the following local weighted Gagliardo-Nirenberg inequality in \cite{dleGNnew}.

\btheo{GNlocalog1m} Suppose LM.1)-LM.2), A.1). Let $U\in C^2(\Og,\mbox{dom}(K))$ and satisfy
\beqno{boundaryzm}\myprod{\og_*\og_0^2 \Fg^2(U)\DKTmU DU,\vec{\nu}}=0\eeq on $\partial\Og$ where $\vec{\nu}$ is the outward normal vector of $\partial\Og$.  Let $\myPi(x):=\LLg^{p+1}(U(x))\Fg^{-p}(U(x))$ and assume that $[\myPi^{\ag}]_{\bg+1}$ is finite for some $\ag>2/(p+2)$ and $\bg<p/(p+2)$.

Then, for any $\eg>0$ there are constants $C, C([\myPi^{\ag}]_{\bg+1})$ such that
\beqno{GNlocog11m}I_{1,*}\le  \eg I_1+\eg^{-1}C\|K(U)\|_{BMO(\mu)}^2[I_2+\myIbar+C([\myPi^{\ag}]_{\bg+1}) [I_2+\myIbar+\breve{I}_{0,*}]].\eeq
Here, $C$ also depends on $C_{P},C_\mu$ and $C_H$.

\etheo

For our purpose in this paper we need only a special case of \reftheo{GNlocalog1m} where $\og$ satisfies AR) so that the Poincar\'e and Hardy inequalities in LM.1) and L.M.2) are verified ($N,s_*=n$). In addition, let $\Og,\Og_*$ be concentric balls $B_s,B_t$, $0<s<t$.
We let $\og_*$ be a cutoff function for $B_s,B_{t}$:  $\og_*$ is a $C^1$ function satisfying $\og_*\equiv1$ in $B_s$ and $\og_*\equiv0$ outside $B_t$ and $|D\og_*|\le 1/(t-s)$. The condition \mref{boundaryzm} of the above theorem is clearly satisfied on the boundary of $\Og=B_t$. We also consider only the case $\Fg(U)\sim|\LLg_U(U)|$.

We then have the following corollary.

\bcoro{GNlocalog1mcoro} Suppose that AR) and A.1) holds for   $\Fg(U)=|\LLg_U(U)|$. Accordingly, define  $\myPi_p(x):=\LLg^{p+1}(U(x))|\LLg_U(U(x))|^{-p}$ and let $B_t(x_0)$ be any ball in $\Og$ and assume that \bdes\item[A.2)]  $[\myPi_p^{\ag}]_{\bg+1,B_t(x_0)}$ is finite for some $\ag>2/(p+2)$ and $\bg<p/(p+2)$.\edes

 We denote (compare with \mref{Idefm}-\mref{I0*m}) \beqno{Idefmt} I_{0}(t,x_0):=\iidmu{B_t(x_0)}{\LLg^2(U)|DU|^{2p}},\; I_1(t,x_0):=\iidmu{B_t(x_0)}{|\LLg_U(U)|^2|DU|^{2p+2}},\;\eeq
\beqno{I0*mt} 
I_2(t,x_0):=\iidmu{B_t(x_0)}{\LLg^2(U)|DU|^{2p-2}|D^2U|^2}.\eeq

Then, for any $\eg>0$ and any ball $B_s(x_0)$, $0<s<t$, there are constants $C, C([\myPi_p^{\ag}]_{\bg+1,B_t(x_0)})$ with $C$ also depending on $C_{PS},C_\mu$ and $C_H$ such that for $$C_{\eg,U,\myPi}=\eg+\eg^{-1}C\|K(U)\|_{BMO(B_t(x_0),\mu)}^2[1+C([\myPi_p^{\ag}]_{\bg+1,B_t(x_0)})]$$ we have
\beqno{GNlocog11mcoro}I_{1}(s,x_0)\le  C_{\eg,U,\myPi}[I_1(t,x_0)+I_2(t,x_0)+(t-s)^{-2}I_{0}(t,x_0)].\eeq 

\ecoro

\brem{GNHrem}We can see that the condition H) implies the condition A.1) in \reftheo{GNlocalog1m}, and then \refcoro{GNlocalog1mcoro}
with $\LLg(u)=\llg^\frac12(u)$ and $\Fg(u)=|\LLg_u(u)|$, \mref{GNlocog11mcoro} is then applicable. Indeed, the assumption \mref{kappamain} in this case is \mref{kappamainz}.  It is not difficult to see that the assumption in f.2) that  $|\llg_{uu}(u)|\llg(u)\lesssim |\llg_u(u)|^2$ and \mref{kappamainz} imply $|\Fg_u(u)||\mathbb{K}(u)|\lesssim \Fg(u)$, which gives \mref{logcondszmain} of A.1). Hence, A.1) holds by H). In particular, if $\llg$ has a polynomial growth in $u$, i.e. $\llg(u)\sim (\llg_0+|u|)^k$ for some $k\ne0$ and $\llg_0\ge0$, then H) reduced to the simple condition $|\mathbb{K}(u)|\lesssim|u|$.
\erem

\section{Proof of The Main Theorem}\eqnoset\label{w12est}

In this section, we prove \reftheo{gentheo1}. We consider the following system
\beqno{gensys}\left\{\barr{l}u_t -\Div(A(x,u)Du)=\hat{f}(x,u,Du),\quad x\in \Og,\\\mbox{$u=0$ or $\frac{\partial u}{\partial \nu}=0$ on $\partial \Og$}. \earr\right.\eeq

We imbed this system in the following family of systems 
\beqno{gensysfam}\left\{\barr{l}u_t -\Div(A(x,\sg u)Du)=\hat{f}(x,\sg u,\sg Du),\quad x\in \Og, \sg\in[0,1],\\\mbox{$u=0$ or $\frac{\partial u}{\partial \nu}=0$ on $\partial \Og$}. \earr\right.\eeq

The proof of \reftheo{gentheo1}, which asserts the existence of strong solutions $u$ to \mref{gensys}, relies on the Leray Schauder fixed point index theorem. Such a strong solution $u$ of \mref{gensys} is a fixed point of a nonlinear map defined on an appropriate Banach space $\mX$. The proof will be based on several lemmas and we will sketch the main steps below.

We will show in \reflemm{dleANSpropt0} that there exist  $p>n/2$ and a constant $M_*$ depending only on the constants in A) and F) such that any strong solution $u$ of \mref{gensysfam} will satisfy \beqno{M*def} \sup_{\tau\in(0,T_0)}\|Du\|_{L^{2p}(\Og,\mu)}\le M_*\; \|u_t\|_{L^{q_0}(Q)}\le M_*.\eeq  We will show  that there are positive constants $\ag,M_0$ such that \beqno{M0def} \|u\|_{C^{\ag,\ag/2}(Q)}\le M_0.\eeq

Following \cite{LSU}, for some $q,r\ge1$ we denote by $\mathbf{V}_{q,r}(Q)$ the Banach space of vector valued functions on $Q$ with finite norm
$$\|u\|_{\mathbf{V}_{q,r}(Q)} =\sup_{t\in(0,T_0)}\|u(\cdot,t)\|_{L^2(\Og)}+\|Du\|_{q,r,Q},$$ where
$$\|v\|_{q,r,Q}:=\left(\int_0^{T_0}\left(\iidx{\Og}{|v(x,t)|^q}\right)^\frac{r}{q}dt\right)^\frac{1}{r}.$$

For $\sg\in[0,1]$ and any $u\in \RR^m$ and $\zeta\in\RR^{mn}$ we  define  the vector valued functions $F^{(\sg)}$ and $f^{(\sg)}$ by \beqno{Bfdef}F^{(\sg)}(x,u,\zeta):=\int_{0}^{1}\partial_\zeta F(\sg,u,t\zeta)\,dt,\quad f^{(\sg)}(x,u):=\int_{0}^{1}\partial_u F(\sg,x,tu,0)\,dt.\eeq

 For   any given $u,w\in \mathbf{V}_{q,r}(Q)$ we write \beqno{Bfdefalin}\mathbf{\hat{f}}(\sg,x, u,w)=F^{(\sg)}(x,u,Du)Dw+f^{(\sg)}(x,u)w+\hat{f}(x,0,0).\eeq

We will define a suitable Banach space $\mX$  and for each $u\in \mX$ we consider the following linear systems, noting that $\mathbf{\hat{f}}(\sg,x,u,w)$ is linear in $w,Dw$
\beqno{Tmapdef}\left\{\barr{l} w_t-\Div(A(x,\sg u)Dw)=\mathbf{\hat{f}}(\sg,x,u,w)\quad (x,t)\in \Og\times(0,T_0), \\\mbox{$w=0$ or $\frac{\partial w}{\partial \nu}=0$ on $\partial \Og\times (0,T_0)$},\\
w(x,0)=U_0(x)\mbox{ on } \Og. \earr\right.\eeq

We will show that the above system has a unique weak solution $w$ if $u$ satisfies \mref{M0def}. We then define $T_\sg(u)=w$ and apply the Leray-Schauder fixed point theorem to establish the existence of a fixed point of $T_1$.
It is clear from \mref{Bfdefalin} that $\hat{f}(x,\sg u,\sg Du)=\mathbf{\hat{f}}(\sg,x,u,u)$.  
Therefore,  from  the definition of $T_\sg$ we see that a fixed point of $T_\sg$ is a weak solution of \mref{gensysfam}. By an appropriate choice of $\mX$, we will show that these fixed points are strong solutions of \mref{gensysfam}, and so a fixed point of $T_1$ is a strong solution of \mref{gensys}.

From the proof of Leray-Schauder fixed point theorem in \cite[Theorem 11.3]{GT}, we need to  find some ball $B_M$ of radius $M$ and centered at $0$ of $\mX$ such that $T_\sg: \bar{B}_M\to \mX$ is compact and  that $T_\sg$ has no fixed point on the boundary of $B_M$.  The topological degree $\mbox{ind}(T_\sg, B_M)$ is then well defined and invariant by homotopy so that $\mbox{ind}(T_1, B_M)=\mbox{ind}(T_0, B_M)$. It is easy to see that the latter is nonzero because the linear system $$\left\{\barr{l}u_t -\Div(A(x,0)Du)=\mathbf{\hat{f}}(x,0,0)\quad x\in \Og\times(0,T_0), \\\mbox{$u=0$ or $\frac{\partial u}{\partial \nu}=0$ on $\partial \Og\times(0,T_0)$, $u(x,0)=U_0(x)$ on $\Og$}, \earr\right.$$ has a unique solution in $B_M$. Hence, $T_1$ has a fixed point in $B_M$.

Therefore, the theorem is proved as we will establish the following claims.

\bdes \item[Claim 1] There exist a Banach space $\mX$ and $M>0$ such that the map $T_\sg:\bar{B}_M\to\mX$ is well defined and compact. 
\item[Claim 2] $T_\sg$ has no fixed point on the boundary of $\bar{B}_M$. That is,  $\|u\|_\mX< M$ for any fixed points of $u=T_\sg(u)$.
\edes

The following lemma defines the space $\mX$, the map $T_\sg$ and establishes Claim 1.

\blemm{Tmaplem} Suppose that there exist $p>n/2$, $q_0>1$ and a constant $M_*$  such that any strong solution $u$ of \mref{gensysfam} satisfies  \beqno{M*defz} \sup_{\tau\in(0,T_0)}\|Du\|_{W^{1,2p}(\Og,\mu)}\le M_*,\; \|u_t\|_{L^{q_0}(Q)}\le M_*.\eeq Then, there exist  $M,\bg>0$ and $q,r\ge1$ such that for $\mX=C^{\bg,\bg/2}(Q,\RR^m)\cap \mathbf{V}_{q,r}(Q)$ the map $T_\sg:\bar{B}_M\to\mX$ is well defined and compact for all $\sg\in[0,1]$. Moreover, $T_\sg$ has no fixed points on $\partial B_M$. \elemm

\bproof 
For some constant $M_0>0$ we consider $u:Q\to\RR^m$ satisfying 
\beqno{Xdefstartzz} \sup_{\tau\in(0,T_0)}\|u\|_{C(\Og)}\le M_0,\; \itQmuz{Q}{|Du|^2}\le M_0,\eeq
and write the system \mref{Tmapdef} as a linear parabolic system for $w$
\beqno{LSUsys} w_t=\Div(\mathbf{a}(u) Dw) + \mathbf{b}(u)Dw+\mathbf{g}(u)w+\mathbf{f},\eeq  where $ \mathbf{a}(x,t)=A(x,\sg u)$, $\mathbf{b}(x,t)=F^{(\sg)}(x,u,Du)$, $\mathbf{g}(x,t)=f^{(\sg)}(x,u)$, and $\mathbf{f}(x)=\hat{f}(x,0,0)$. 
The matrix $\mathbf{a}(u)$ being regular elliptic with uniform ellipticity constants by A), AR) if $u$ is bounded. We recall the following well known result in \cite[Chapter VII]{LSU}. If there exist positive constants $m$ and $q,r$ such that (see the condition (1.5) in \cite[Chapter VII]{LSU}) \beqno{LSUcond}\|\mathbf{b}(u)\|_{q,r,Q},\; \|\mathbf{g}(u)\|_{q,r,Q},\; \|\mathbf{f}\|_{q,r,Q}\le m,\; \mbox{$1/r+n/(2q)=1$, $q\ge n/2$ and $r\ge1$},\eeq then the system \mref{LSUsys} satisfies the assumptions of Theorem 1.1 in \cite[Chapter VII]{LSU} which asserts that \mref{Tmapdef} has a unique weak solution $w$.

 Moreover, as the initial condition $w(\cdot,0)=U_0(x)$ belongs to $W^{1,p_0}(\Og)$ and then $C^{\bg_0}(\Og)$ for $\bg_0=1-n/p_0>0$, a combination of Theorems 2.1 and 3.1 in \cite[Chapter VII]{LSU} shows that $w$ belongs to $C^{\ag_0,\ag_0/2}(\bar{Q},\RR^m)$ for some $\ag_0>0$ depending only on $\bg_0$,$\|u\|_\infty$ and $m$.

Next, we will show that \mref{LSUcond} holds  by F) and \mref{Xdefstartzz}. We consider the two cases  f.1) and  f.2). If f.1) holds then from the definition \mref{Bfdef} there is a constant $C(|u|)$ such that  $$|\mathbf{b}(x,t)|=|F^{(\sg)}(x,u,\zeta)|\le C(|u|),\; |\mathbf{g}(x,t)|=|f^{(\sg)}(x,u)|\le C(|u|).$$
 
 From  \mref{Xdefstartzz}, we see that $\sup_{\tau\in(0,T_0)}\|u\|_\infty\le M_0$ and so there is  a constant $m$ depending on $M_0$ such that \mref{LSUcond} holds for any $q$ and $n$.

If f.2) holds then
\beqno{Bfdeff2}|F^{(\sg)}(x,u,\zeta)|\le C(|u|)|\zeta|,\; |f^{(\sg)}(x,u)|\le C(|u|).\eeq Therefore, $\|\mathbf{b}\|_{L^2(Q)}$ is bounded by $C\|Du\|_{L^2(Q)}$. Again, if  $n\le3$ then  \mref{Xdefstartzz} implies the condition \mref{LSUcond} for $q=2$.

In both cases,  \mref{LSUsys} (or \mref{Tmapdef}) has a unique weak solution $w$. 
We then define $T_\sg(u)=w$. Moreover, as we explained earlier, $w\in C^{\ag_0,\ag_0/2}(Q)$ for some $\ag_0>0$ depending on $M_0$. 

We now consider a fixed point $u$ of $T_\sg$. By \reflemm{Dulocbound} following this proof we see that $u$ is a strong solution and we can use the assumption \mref{M*defz}.
The first bound in the assumption \mref{M*defz} implies $u$ is H\"older continuous in $x$. This and the integrability of $u_t$ in the second bound of the assumption and \cite[Lemma 4]{NecasSverak} provide positive constants $\ag,M_1$ such that any strong solution $u$ of \mref{gensysfam} satisfies $\|u\|_{C^{\ag,\ag/2}(\Og)}\le M_1$.  Also, the assumption AR) implies that $\llg(u),\og$ are bounded from below, yield that $\|Du\|_{L^2(Q)}\le C(C_0)$. Thus, there is a constant $M_1$, depending on $M_*,C_0$ such that any strong solution $u$ of \mref{gensysfam} satisfies 
\beqno{Xdefstart} \|u\|_{C^{\ag,\ag/2}(Q)}\le M_1,\; \|Du\|_{L^2(Q)}\le M_1.\eeq

It is well known that there is a constant $c_0>1$, depending on $\ag,T_0$ and the diameter of $\Og$, such that $\|\cdot\|_{C^{\bg,\bg/2}(Q)}\le c_0\|\cdot\|_{C^{\ag,\ag/2}(Q)}$ for all $\bg\in(0,\ag)$. We now let $M_0$, the constant in \mref{Xdefstartzz}, be $M=(c_0+1)M_1$.

Define $\mX=C^{\bg,\bg/2}(Q)\cap V^{1,0}(Q)$ for some positive $\bg<\min\{\ag,\ag_0\}$, where $\ag_0$ is the H\"older continuity exponent for solutions of \mref{LSUsys}, and $$V^{1,0}(Q):=\{u\,:\, Du \in L^2(Q)\}.$$

The space $\mX$ is equipped with the norm $\|u\|_\mX = \max\{\|u\|_{C^{\bg,\bg/2}(Q)},\|Du\|_{L^2(Q)}\}$ and consider the ball $B_M$ in $\mX$ centered at $0$ with radius $M$.

We now see that $T_\sg$ is well defined and maps the ball $\bar{B}_M$ of $\mX$ into $\mX$. Moreover, from the definition $M=(c_0+1)M_1$, it is clear that $T_\sg$ has no fixed point on the boundary of $B_M$ because such a fixed points $u$ satisfies \mref{Xdefstart} which implies $\|u\|_{\mX}\le c_0M_1<M$. 
 
Finally, we need only show that $T_\sg$ is compact.  If $u$ belongs to a bounded set $K$ of $\bar{B}_M$ then $\|u\|_\mX\le C(K)$ for some constant $C(K)$ and there is a constant $C_1(K)$  such that $\|T_\sg(u)\|_{C^{\ag_0,\ag_0/2}(Q)}=\|w\|_{C^{\ag_0,\ag_0/2}(Q)}\le C_1(K)$. Thus $T_\sg(K)$ is compact in $C^{\bg,\bg/2}(Q)$ because $\bg<\ag_0$.  So, we need only show that $T(K)$ is precompact in $V^{1,0}(Q)$.  We will discuss only the quadratic growth case where \mref{Bfdeff2} holds because the case $\hat{f}$ has linear growth is similar and easier.

First of all, for $u\in K$ we easily see that $\|Dw\|_{L^2(Q)}$ is uniformly bounded by a constant depending on $K$. The argument is standard by testing the linear system \mref{Tmapdef} by $w$ and using the boundedness of $\|w\|_{L^\infty}$ and $\|u\|_{L^\infty}$, \mref{Bfdeff2}, AR) and Young's inequality.

Let $\{u_n\}$ be a sequence in $K$ and $w_n=T_\sg(u_n)$.  We have, writing $W=w_n-w_m$  $$W_t-\Div(A(x,\sg u_n)DW)=\Div(\ag_{m,n}Dw_m)+\Psi_{m,n},$$ where $\ag_{m,n}=(A(x,\sg u_n)-A(x,\sg u_m)$ and $\Psi_{m,n}$ is defined by $$F^{(\sg)}(x,u_n,Du_m)Dw_n-F^{(\sg)}(x,u_m,Du_m)Dw_m+f^{(\sg)}(x,u_n)u_n-f^{(\sg)}(x,u_m)u_m.$$

Testing the above system with $W$ and using AR) and the fact that $W(x,0)=0$, we have for $dz=dxdt$
$$\llg_*(K)\mu_*\itQ{Q}{|DW|^2} \le \itQ{Q}{[|\ag_{m,n}||Dw_m||DW|+|\Psi_{m,n}||W|]}.$$
By Young's inequality, we find a constant $C$ depending on $K$ and $\mu_*$ such that
$$\itQ{Q}{|DW|^2}\le C\itQ{Q}{[(|\ag_{m,n}||Dw_m|)^2} +\sup_Q|W|\|\Psi_{m,n}\|_{L^1(Q)}.$$
By \mref{Bfdeff2}, it is clear that   $|\Psi_{n,m}|\le C(K)[(|Du_n|+|Du_m|)(|Dw_n|+|Dw_m|)+1]$.
Using the fact that $\|Dw_n\|_{L^2(Q)}$ and $\|Du_n\|_{L^2(Q)}$ are uniformly bounded, we see that $\|\Psi_{m,n}\|_{L^1(Q)}$ is bounded. Hence, $$\itQ{Q}{|Dw_n-Dw_m|^2}\le C(K)\max\{\sup_\Og|A(x,\sg u_n)-A(x,\sg u_m)|,\sup_\Og|w_n-w_m|\}.$$

Since $u_n,w_n$ are bounded in $C^{\bg,\bg/2}(Q))$, passing to subsequences we can assume that $u_n,w_n$ converge in $C^0(Q)$. Thus, $\|A(x,\sg u_n)-A(x,\sg u_m)\|_\infty,\|w_n-w_m\|_\infty\to0$. We then see from the above estimate that $Dw_n$ converges in $L^{2}(\Og)$. Thus, $T_\sg(K)$ is precompact in $V^{1,0}(Q)$.

Hence, 
$T_\sg:\mX\to\mX$ is a compact map. The proof is complete. \eproof

We now turn to Claim 2, the hardest part of the proof, and provide a uniform estimate for the fixed points of $T_\sg$ and justify the key assumption \mref{M*defz} of \reflemm{Tmaplem}. The proof is complicated and will be devided into many lemmas described as follows.
\begin{itemize} \item  \reflemm{Dulocbound} is quite standard and shows that the fixed points of $T_\sg$ are strong solutions.
\item \reflemm{dleANSenergy} follows \cite[Lemma 3.2]{dleANS} and establishes an energy estimate of $Du$. In \reflemm{dleANSprop}, the assumptions H) and M.1) then allow us to apply the local Gagliardo-Nirenberg inequality \mref{GNlocog11mcoro}  to obtain a better estimate.
\item \reflemm{dleANSpropt0} and \reflemm{selfimprovlem} then show that the estimate in \reflemm{dleANSprop} is self-improving to obtain the key estimate \mref{M*defz}.
\end{itemize}

Hence, we first have the following lemma. 

\blemm{Dulocbound} A fixed point of $T_\sg$ is also a strong solution of \mref{gensysfam}. 
 \elemm
\bproof If $u$ is a fixed point of $T_\sg$ in $\mX$ then it solves \mref{gensysfam} weakly and is continuous. Thus, $u$ is bounded and belongs to $VMO(Q)$. By AR), the system \mref{gensysfam} is regular elliptic. We can adapt the proof in \cite{GiaS}.
If $\hat{f}$ satisfies a quadratic growth in $Du$ then,  because $u$ is bounded, the condition \cite[(0.4)]{GiaS} that $|\hat{f}|\le a|Du|^2+b$ is satisfied here. The proof of \cite[Theorems 2.1 and 3.2]{GiaS} assumed the 'smallness condition' (see \cite[(0.6)]{GiaS}) $2aM<\llg_0$, where $M=\sup |u|$. This 'smallness condition' was needed because only weak bounded solutions, which are not necessarily continuous, were considered in \cite{GiaS}. In our case,  $u$ is continuous so that we do not require this 'smallness condition'.  Indeed, a careful checking of the arguments of the proof in \cite[Lemma 2.1 and page 445]{GiaS} shows that if $R$ is small and one knows that the solution $u$ is continuous then these argument still hold as long we can absorb the integrals involving $|Du|^2, |Dw|^2$ (see the estimate after \cite[(3.7)]{GiaS}) on the right hand sides to the left right hand sides of the estimates. Thus, \cite[Theorems 2.1 and 3.2]{GiaS} apply to our case and yield that $u\in C^{a,a/2}(Q)$ for all $a\in(0,1)$ and that, since $A(x,u)$ is differentiable,  $Du$ is locally H\"older continuous in $Q$. Therefore, $u$ is also a strong solution. \eproof

Thanks to \reflemm{Dulocbound}, we need only consider a strong solution  $u$ of \mref{gensysfam} and establish \mref{M*defz} for some $p>n/2$. Because the data of \mref{Tmapdef} satisfy the structural conditions A), F) with the same set of constants and the assumptions of the theorem are assumed to be uniform for all $\sg\in[0,1]$, we will only present the proof for the case $\sg=1$ in the sequel.

Let $u$ be a strong solution  of \mref{e1} on $\Og$.  We begin with an energy estimate for $Du$. For $p\ge1$ and any ball $B_s$ with center $x_0\in\bar{\Og}$ we denote $\Og_s=B_s\cap\Og$, $Q_s=\Og_s\times(0,T_0)$ and 

\beqno{Adef}\ccA_p(s)=\sup_{\tau\in (0,T_0)}\iidx{\Og_s\times {\tau}}{|Du|^{2p}},\eeq
\beqno{Hdef}\ccH_{p}(s):=
\itQmuz{Q_s}{\llg(u)|Du|^{2p-2}|D^2u|^2},\eeq
\beqno{Bdef}\ccB_{p}(s):= \itQmuz{Q_s}{\frac{|\llg_u(u)|^2}{\llg(u)}|Du|^{2p+2}},\eeq \beqno{Cdef}\ccC_{p}(s):=\itQmuz{Q_s}{(|f_u(u)|+\llg(u))|Du|^{2p}},\eeq and \beqno{cIdef}\mccIee_{\og,p}(s):=\itQmuz{Q_s}{(\llg(u)|Du|^{2p}|D\og_0|^2+|f(u)||Du|^{2p-1}|D\og_0|\og_0)}.\eeq

The following lemma establishes an energy estimate for $Du$.

\blemm{dleANSenergy}  Assume A), F). 
Let $u$ be any strong solution  of \mref{gensys} on $\Og$ and  $p$ be any number in $[1,\max\{1,n/2\}]$.

There is a constant $C$, which depends only on the parameters in A) and F), such that for any two concentric balls $B_s,B_t$ with center $x_0\in\bar{\Og}$ and $s<t$ \beqno{keydupANSenergy}\ccA_p(s)+\ccH_{p}(s)\le C\ccB_{p}(t)+C(1+(t-s)^{-2})[\ccC_{p}(t)+\mccIee_{\og,p}(t)]+\||DU_0|^{2p}\|_{L^1(\Og_t)}.\eeq

\elemm

\bproof The proof is similar to the energy estimate of $Du$ for the parabolic case in \cite[Lemma 3.2]{dleANS}. Roughly speaking, we differentiated the system in $x$ to obtain
\beqno{ga2zzz} (Du)_t-\Div(A(x,u)D^2u+A_u(x,u)DuDu+A_x(x,u)Du)=D\hat{f}(x,u,Du).\eeq

For any two concentric balls $B_s,B_t$, with $s<t$, let $\psi$ be a cutoff function for $B_s,B_{t}$. That is, $\psi$ is a $C^1$ function satisfying $\psi\equiv1$ in $B_s$ and $\psi\equiv0$ outside $B_t$ and $|D\psi|\le 1/(t-s)$. Consider any given triple $t_0,T,T'$ satisfying $0<t_0<T<T'\le T_0$ and $\eta$  being a cutoff function for $(T-t_0,T'),(T,T')$. We then
test \mref{ga2zzz} with $|Du|^{2p-2}Du\psi^2\eta$ and obtain, using integration by parts and Young's inequality

\beqno{keyenergystart}\barr{ll}\lefteqn{\sup_{t\in(T,T')}\iidx{\Og_s}{|Du|^{2p}}+
	\dspl{\int}_{T-t_0}^{T'}\iidx{\Og_t}{\llg(u)|Du|^{2p-2}|D^2u|^2\psi^2\eta\og}d\tau\le }\hspace{.2cm}&\\& C\dspl{\int}_{T-t_0}^{T'}\iidx{\Og_t}{[\frac{|\llg_u(u)|^2}{\llg(u)}|Du|^{2p+2}+|D\psi|^{2}\llg(u)|Du|^{2p}]\eta\og} d\tau\\& +C\dspl{\int}_{T-t_0}^{T'}\iidx{\Og_t}{[|A_x(x,u)||Du|^{2p-1}|D^2u|\psi^2+|D\hat{f}(x,u,Du)||Du|^{2p-1}\psi^2]}d\tau\\&+Ct_0^{-1}\dspl{\int}_{T-t_0}^{T}\iidx{\Og_t}{|Du|^{2p}}d\tau.\earr\eeq

Here, integrals in the first line of \mref{keyenergystart} result from the same argument in the proof of \cite[Lemma 3.2]{dleANS} using the spectral gap condition SG) we are assuming here (see also \cite[Lemma 6.5]{dleAMJ}).  The integrals in the second and third lines can be estimated by simple uses of Young's inequality and the condition F) as in \cite{dleANS,dleAMJ}. Finally, We formally let $T,t_0\to0$ in the last integral, which will be justified below, to obtain \mref{keydupANSenergy}.

Using the difference quotience operator $\dg_h$ instead of $D$ in \mref{ga2zzz}, we obtain
\beqno{ga2zzz1} (\dg_hu)_t=\Div(A(x,u)D(\dg_h u)+\dg_h(A(x,u))Du)+\dg_h\hat{f}(x,u,Du).\eeq 

We test this with $|\dg_hu|^{2p-2}\dg_hu\psi^2\eta$ to obtain a similar version of \mref{keyenergystart} with the operator $D$ being replaced by $\dg_h$. We can integrate the result over $(0,T_0)$ and obtain
$$\barr{ll}\lefteqn{\sup_{t\in(0,T_0)}\iidx{\Og_s}{|\dg_hu|^{2p}}+
	\itQ{Q_{s}}{\llg(u)|\dg_hu|^{2p-2}|D\dg_hu|^2}\le}\hspace{.5cm}&\\& C\itQ{Q_{t}}{[\frac{|\llg_u(u)|^2}{\llg(u)}|Du|^2|\dg_hu|^{2p}+|D\psi|^2\llg(u)|\dg_hu|^{2p}]}+\cdots+C\iidx{\Og_t}{|\dg_hu(x,0)|^{2p}}.\earr$$

Since $u\in C([0,T'),L^{2p}(\Og))$,  we can let $h$ tend to 0 and obtain a similar energy estimate \mref{keydupANSenergy} for $Du$ with $T=t_0=0$ and $\eta\equiv1$. We complete the proof. \eproof

Next, under the condition AR), the density $\og$  supports the Poincar\'e-Sobolev inequality   with  $\pi_*=2n/(n-2)$. By \refrem{GNHrem}, we can apply the local Gagliardo-Nirenberg inequality \mref{GNlocog11mcoro} here. 
Thus, if the condition \mref{Keymu0} of M.1) holds then we combine the energy estimate and \mref{GNlocog11mcoro} to have the following stronger estimate.

\blemm{dleANSprop}  In addition to the assumptions of \reflemm{dleANSenergy}, we suppose that H) and M.1) hold for some $p$. That is, for any given  
$\mu_0>0$ there exist a constant $C_0$ and a positive $R_{\mu_0}$ sufficiently small in terms of the constants in A) and F) such that
\beqno{Keymu01} \sup_{x_0\in\bar{\Og},\tau\in(0,T_0)}[\myPi_p^{\ag}]_{\bg+1,\Og_R(x_0)}\le C_0,\;\sup_{x_0\in\bar{\Og},\tau\in(0,T_0)}\|K(u)\|_{BMO(\Og_{R}(x_0),\mu)}^2 \le \mu_0.\eeq

Then for sufficiently small $\mu_0$ there is a constant $C$ depending only on the parameters of A) and F) such that for $2R<R_{\mu_0}$ we have \beqno{keydupANSppb}\ccA_p(R)+\ccB_{p}(R)+\ccH_{p}(R)\le C(1+R^{-2})[\ccC_{p}(2R)+\mccIee_{\og,p}(2R)] +\||DU_0|^{2p}\|_{L^1(\Og_{2R})}.\eeq

 \elemm

\bproof  Recall the energy estimate \mref{keydupANSenergy} in \reflemm{dleANSenergy} 
\beqno{keyiteration1} \ccA_p(s)+\ccH_{p}(s)\le C\ccB_{p}(t)+C(1+(t-s)^{-2})[\ccC_{p}(t)+\mccIee_{\og,p}(t)]+\||DU_0|^{2p}\|_{L^1(\Og_{t})}, \; 0<s<t.\eeq 

We apply \refcoro{GNlocalog1mcoro} to estimate $\ccB_{p}(t)$, the integral  on the right hand side of \mref{keyiteration1}. We let $\LLg(u)=\llg^\frac{1}{2}(u)$ in \refcoro{GNlocalog1mcoro} and note that $\myPi_p$ defined there is now comparable to the $\myPi_p=\llg^{p+\frac12}(u)|\llg_u(u)|^{-p}$ in M.1). We compare the definitions \mref{Idefmt} and \mref{I0*mt} with those in \mref{Hdef}-\mref{Cdef} to see that for $U(x)=u(x,\tau)$ with $\tau\in(0,T_0)$ $$\ccB_{p}(t)=\int_0^{T_0}I_1(t,x_0)d\tau,\;\ccC_{p}(t)=\int_0^{T_0}I_0(t,x_0)d\tau,\;\ccH_{p}(t)=\int_0^{T_0}I_2(t,x_0)d\tau.$$
Hence, for any $\eg>0$ we can use \mref{GNlocog11mcoro} obtain a constant $C$ such that (using the bound $[\myPi_p^{\ag}]_{\bg+1,B_{R_{\mu_0}}(x_0)\cap\Og}\le C_0$ and the definitions of $\mu_0$ in \mref{Keymu01}  and $C(\eg,U,\myPi)$  in \refcoro{GNlocalog1mcoro})
$$I_{1}(s,x_0)\le  C_{\eg,U,\myPi}[I_1(t,x_0)+I_2(t,x_0)+(t-s)^{-2}I_{0}(t,x_0)].$$
Integrating the above over $(0,T_0)$ to get
$$ \ccB_{p}(s) \le \eg\ccB_{p}(t)+C\eg^{-1}\mu_0\ccH_{p}(t)+C\eg^{-1}\mu_0(t-s)^{-2}\ccC_{p}(t)\quad 0<s<t\le R_{\mu_0}.$$

Define $F(t):=\ccB_{p}(t)$, $G(t):=\ccH_{p}(t)$, $g(t):=\ccC_{p}(t)$ and $\eg_0=\eg+C\eg^{-1}\mu_0$. The above yields \beqno{keyiteration2}F(s)\le \eg_0[F(t)+G(t)] +C(t-s)^{-2}g(t).\eeq 

Now, for $h(t):=\mccIee_{\og,p}(t)+\||DU_0|^{2p}\|_{L^1(\Og_t)}$ the energy estimate \mref{keyiteration1}  implies
 \beqno{Giter}G(s)\le C[F(t)+(1+(t-s)^{-2})(g(t)+h(t))].\eeq

As $\eg_0=\eg+C\eg^{-1}\mu_0$, it is clear that we can choose and fix some $\eg$ sufficiently small and then $\mu_0$ small in terms of $C,\eg$ to have $2C\eg_0<1$. Thus, if $\mu_0$ is sufficiently small in terms of the  constants in A),F), then we can apply a simple iteration argument \cite[Lemma 3.11]{dleANS} to the two inequalities \mref{keyiteration2} and \mref{Giter} and obtain for $0<s<t\le R_{\mu_0}$ $$ F(s)+G(s) \le C(1+(t-s)^{-2})[g(t)+h(t)].$$ 

For any $R<R_{{\mu_0}}/2$ we take $t=2R$ and $s=\frac32R$ in the above to obtain $$\ccB_{p}(\frac32R)+\ccH_{p}(\frac32R)\le C(1+R^{-2})[\ccC_{p}(2R)+\mccIee_{\og,p}(2R)+\||DU_0|^{2p}\|_{L^1(\Og_t)}].$$ Combining this and \mref{keyiteration1} with $s=R$ and $t=\frac32R$, we see that
$$\ccA_p(R)+\ccB_{p}(R)+\ccH_{p}(R)\le C(1+R^{-2})[\ccC_{p}(2R)+\mccIee_{\og,p}(2R)+\||DU_0|^{2p}\|_{L^1(\Og_t)}].$$
This is \mref{keydupANSppb} and the proof is complete. \eproof

Finally, we have the following lemma giving a uniform bound for strong solutions.

\blemm{dleANSpropt0}  Assume as in \reflemm{dleANSprop} and AR). We assume also the integrability condition M.0).
Then there exist $p>n/2$, $q_0>1$ and a constant $M_*$ depending only on the parameters of A) and F), $\mu_0$, $R_{\mu_0}$, $C_0$ and the geometry of $\Og$  such that 
\beqno{keydupANSt0}\sup_{\tau\in(0,T_0)}\iidmu{\Og}{|Du(\cdot,\tau)|^{2p}}\le   M_*,\eeq
\beqno{keydutANSt0}\|u_t\|_{L^{q_0}(Q)}\le   M_*.\eeq
\elemm

\bproof First of all, by the condition AR), there is a constant $C_\og$ such that $|D\og_0|\le C_\og\og_0$ and therefore we have from the the definition \mref{cIdef} that $$\mccIee_{\og,p}(s)\le C_\og\itQmuz{Q_s}{(\llg(u)|Du|^{2p}+f(u)|Du|^{2p-1})\og_0^2}.$$ By Young's inequality, $f(u)|Du|^{2p-1}\lesssim |f_u(u)||Du|^{2p}+(f(u)|f_u(u)|^{-1})^{2p}|f_u(u)|$. It follows from the assumption  \mref{specfcond} that $(f(u)|f_u(u)|^{-1})^{2p}|f_u(u)|\lesssim (|u|+1)^{2p}|f_u(u)|$. We then have from \mref{keydupANSppb} that
\beqno{keydupANSppbbb}\ccA_p(R)+\ccB_{p}(R)+\ccH_{p}(R)\le C(1+R^{-2})[\ccC_{p}(2R)+\mccIee_{*,\og,p}(2R)+C_0],\eeq \beqno{cIdef1}\mccIee_{*,\og,p}(s):=\itQmuz{Q_s}{|u|^{2p}|f_u(u)|}.\eeq

The main idea of the proof is to show that \mref{keydupANSppbbb} is self-improving in the sense that if it is true for some exponent $p\ge1$ then it is also true for $\cg_*p$ with some fixed $\cg_*>1$ and $R$ being replaced by $R/2$. To this end, assume that for some $p\ge1$ we can find a constant $C(C_0,R,p)$ such that \beqno{piter1} \ccC_{p}(2R)+\mccIee_{*,\og,p}(2R)\le C(C_0,R,p),\eeq which and \mref{keydupANSppbbb} and the definitions of $\ccB_{p}(R),\ccH_{p}(R),\ccC_p(R)$ yield that $$ \ccA_p(R)+ \itQmuz{Q_R}{[\llg(u)|Du|^{2p}+\llg(u)|Du|^{2p-2}|D^2u|^2+\Fg^2(u)|Du|^{2p+2}]}\le C(C_0,R,p),$$
where $\Fg(u)=|(\llg^\frac12(u))_u|$.  The above two estimates  yield for $V=\llg^\frac12(u)|Du|^{p}$ \beqno{keydupANSpp1} \sup_{\tau\in(0,T_0)}\iidx{\Og_R}{|Du|^{2p}}+
\itQmuz{Q_R}{[V^2+|DV|^2]}\le C(C_0,R,p).\eeq

In the technical \reflemm{selfimprovlem} following this proof, we will show that if M.0) and \mref{piter1} hold for some $p\ge1$ then together with its consequence \mref{keydupANSpp1} provide some $\cg_*>1$ such that  \mref{piter1} holds again for the new exponent $\cg_{*} p$  and $R/2$. 

By the assumption \mref{hypoiterpis1}, \mref{piter1} holds for $p=1$. It is now clear that, as long as the energy estimate \mref{keydupANSenergy} is valid by \reflemm{dleANSenergy}), we can repeat the argument $k_0$ times to find  a number $p>n/2$ such that \mref{piter1} and then its consequence \mref{keydupANSpp1} hold. It follows  that 
there is a constant $C$ depending only on the parameters of A) and F), $\mu_0$, $R_{\mu_0}$ and $k_0$ 
such that for some $p>n/2$ we obtain from \mref{keydupANSpp1} that
\beqno{keydupANSt0z}\sup_{(0,T_0)}\iidmu{\Og_{R_0}}{|Du|^{2p}}\le   C\mbox{ for  $R_0=2^{-k_0}R_{\mu_0}$}.\eeq

Summing the above inequalities over a finite covering of balls $B_{R_0}$ for $\Og$, we find a constant $C$, depending also on the geometry of $\Og$, and obtain the desired estimate \mref{keydupANSt0}. 

Similarly, we obtain from \mref{keydupANSpp1} with $p=1$ that \beqno{keyD2u}\itQmuz{Q}{\llg(u)|D^2u|^2} \le C.\eeq As $u$ is a strong solution, we have $|u_t|\le |\Div(A(x,u)Du)|+|\hat{f}|$ a.e. in $Q$. Therefore, $$\|u_t\|_{L^{q_0}(Q)}\lesssim \|\llg(u)|D^2u|\|_{L^{q_0}(Q)}+\|\llg_u(u)|Du|^2\|_{L^{q_0}(Q)}+\|f(u)\|_{L^{q_0}(Q)}.$$ If $q_0\in(1,2)$ then the first and third norms on the right can be treated by H\"older's inequality and \mref{keyD2u} and the boundedness of $u$, thanks to \mref{keydupANSt0}. For $q_0=p>n/2\ge1$, the second norm is also bounded by \mref{keydupANSt0}. Thus, there is $q_0>1$ such that \mref{keydutANSt0} holds.

The lemma is proved. \eproof

Thus, we need to show that \mref{piter1} is self improving in the following lemma. 

\blemm{selfimprovlem} Assume as in \reflemm{dleANSpropt0}. Suppose that for some $p\ge1$ we can find a constant $C(C_0,R,p)$ such that \beqno{piter1a}  \ccC_{p}(2R)+\mccIee_{*,\og,p}(2R)\le C(C_0,R,p),\eeq then there exists a fixed $\cg_*>1$ such that \beqno{piter1azzz} \ccC_{\cg_*p}(R) +\mccIee_{*,\og,\cg_*p}(R)\le C(C_0,R,p).\eeq
\elemm

In the sequel,  we will repeatedly make use of the following parabolic Sobolev inequality \beqno{paraPS}\itQmuz{Q_R}{v^{2q_{*}}|V|^{2}}\lesssim \sup_I\mypar{\iidmu{\Og_R}{v^2}}^{q_*}\itQmuz{Q_R}{[|DV|^{2}+V^2]}, \mbox{ $q_*=1-\frac2{\pi_*}$}.\eeq 

To see this, we recall the inequality $$\mypar{\iidmu{\Og}{|V|^{\pi_*}}}^\frac1{\pi_*}\lesssim \mypar{\iidmu{\Og}{|DV|^{2}}}^\frac1{2} + \mypar{\iidmu{\Og}{|V|^2}}^\frac{1}{2}$$ which is just a simple consequence of the Poincar\'e-Sobolev inequality PS). For $q_{*}=(1-\frac2{\pi_*})$ we  use H\"older's inequality and the above inequality to have $$\barr{lll}\itQmuz{\Og\times I}{v^{2q_{*}}|V|^{2}}&\le& \dspl{\int}_I\mypar{\iidmu{\Og}{v^2}}^{1-\frac{2}{\pi_*}}\mypar{\iidmu{\Og}{|V|^{\pi_*}}}^\frac{2}{\pi_*}d\tau\\ &\lesssim& \sup_I\mypar{\iidmu{\Og}{v^2}}^{q_*}\dspl{\int_I}\mypar{\iidmu{\Og}{|DV|^{2}}+\iidmu{\Og}{|V|^2}}d\tau.\earr$$
This is \mref{paraPS}.

{\bf Proof of \reflemm{selfimprovlem}:} We recall the integrability condition M.0). Namely, there exists $C_0$ and $r_0>1,\bg_0\in(0,1)$ such that\beqno{llgmainhypb} \sup_{\tau\in(0,T_0)}\||f_u(u)|\llg^{-1}(u)\|_{L^{r_0}(\Og,\mu)},\; \sup_{\tau\in(0,T_0)}\|u^{\bg_0}\|_{L^{1}(\Og,\mu)}\le C_0 ,\eeq
\beqno{hypoiterpis1b}\itQmuz{Q}{(|f_u(u)|+\llg(u))(|Du|^{2}+|u|^2)}\le C_0,\eeq

 We established in the proof of \reflemm{dleANSpropt0} that for $V=\llg^\frac12(u)|Du|^{p}$ \mref{piter1a}   yields  \mref{keydupANSppbbb}, which and the fact that $\og$ is bounded from above imply \beqno{keydupANSpp1a} \sup_{\tau\in(0,T_0)}\iidmu{\Og_R}{|Du|^{2p}}+
\itQ{Q_R}{[V^2+|DV|^2]}\le C(C_0,R,p).\eeq

Let $\cg_{1}=1+q_*$. We write $\llg(u)|Du|^{\cg_12p}=v^{2q_*}V^2$  with $v=|Du|^p$, $V=\llg^\frac12(u)|Du|^{p}$  and apply \mref{paraPS} to get
$$\itQmuz{Q_R}{\llg(u)|Du|^{\cg_{1}2p}}\lesssim \sup_{(0,T_0)}\mypar{\iidmu{\Og_R}{|Du|^{2p}}}^{q_*} \itQmuz{Q_R}{[|DV|^{2}+V^2]}.$$  Therefore, \mref{keydupANSpp1a} implies $$\itQmuz{Q_R}{\llg(u)|Du|^{\cg_{1}2p}}\le C(C_0,R,p).$$

Similarly, we write $|f_u(u)||Du|^{\cg_22p}=v^{2q_*}V^2$ with $v=(|f_u(u)|\llg^{-1}(u)|Du|^{2p(\cg_2-1)})^\frac{1}{2q_*}$ and $V=\llg^\frac12(u)|Du|^{p}$. In order to apply \mref{paraPS} here, we need to estimate the integral of $v^2$ over $\Og_R$. Assuming $\cg_2\in(1,q_*)$ and using H\"older's inequality with the exponent $q_1=\frac{q_*}{q_*-\cg_2+1}$, the integral of $v^2=(|f_u(u)|\llg^{-1}(u)|Du|^{2p(\cg_2-1)})^\frac{1}{q_*}$ is bounded by  $$ \mypar{\iidmu{\Og_R}{(|f_u(u)|\llg^{-1}(u))^{q_1}}}^\frac{1}{q_1}\mypar{\iidmu{\Og_R}{|Du|^{2p}}}^\frac{1}{q_1'}$$ We can find $\cg_2$ close to 1 such that $q_1\le r_0$, which is greater than 1, so that the first integral is bounded by the assumption \mref{llgmainhypb}. The second integral is bounded because of \mref{keydupANSpp1a}.

We now turn to $\mccIee_{*,\og,p}(s)$ defined by \mref{cIdef1} and write $$I_p(s):=\mccIee_{*,\og,p}(s)=\itQmuz{\Og_s}{|f_u(u)||u|^{2p}},\; J_p(s)=\itQmuz{\Og_s}{\llg(u)|u|^{2p}}.$$ 
We will prove that $I_p(R), J_p(R)$
are self improving. We assume first that
\beqno{IJp} I_p(R),\; J_p(R)\le C(C_0,R,p).\eeq

The argument is very similar to the above treatment of the integral of $|f_u(u)||Du|^{\cg_22p}$ with $Du$ being replaced by $|u|$. In fact, the proof for $I_p,J_p$ are almost identical so that we will denote $g(u)=|f_u(u)|$ and consider $I_p$ first. We write $g(u)|u|^{\cg_22p}=v^{2q_*}V^2$ with $v=(g(u)\llg^{-1}(u)|u|^{2p(\cg_2-1)})^\frac{1}{2q_*}$ and $V=\llg^\frac12(u)|u|^{p}$.
We use \mref{paraPS}  to have \beqno{Ipest}I_p(R)\lesssim \sup_I\mypar{\iidmu{\Og_R}{v^2}}^{1-\frac2{\pi_*}}\itQmuz{Q_R}{[|DV|^2+V^2]}.\eeq 

The integral of $v^2=(g(u)\llg^{-1}(u)|u|^{2p(\cg_2-1)})^\frac{1}{q_*}$ over $\Og_R$ is estimated by H\"older's inequality  as before by $$ \mypar{\iidmu{\Og_R}{(g(u)\llg^{-1}(u))^{q_1}}}^\frac{1}{q_1}\mypar{\iidmu{\Og_R}{|u|^{2p}}}^\frac{1}{q_2}$$ Again, the first integral is bounded by \mref{llgmainhypb} as $g(u)=|f_u(u)|$. We consider the second integral and use Sobolev's inequality to have $$\iidmu{\Og_R}{|u|^{2p}}\lesssim \iidmu{\Og_R}{|D(|u|^p)|^2}+\mypar{\iidmu{\Og_R}{|u|^{p\bg}}}^\frac{2}{\bg}.$$ Because $|D(|u|^p)|^2\sim |u|^{2p-2}|Du|^2\le \eg|u|^{2p}+C(\eg)|Du|^{2p}$, we conclude that $$\iidmu{\Og_R}{|u|^{2p}}\lesssim \iidmu{\Og_R}{|Du|^{2p}}+\mypar{\iidmu{\Og_R}{|u|^{p\bg}}}^\frac{2}{\bg}.$$ The first integral on the righ hand side is bounded by \mref{keydupANSpp1a}. Taking $\bg=\bg_0$, the second integral is bounded by the assumption \mref{llgmainhypb}.

Finally, for the last integral in \mref{Ipest} with $V=\llg^\frac12(u)|u|^{p}$ we use the fact that $|\llg_u(u)||u|\lesssim \llg(u)$ and Young's inequality to see that $$\barr{lll}|DV|^2&\lesssim& \llg(u)|D(|u|^p)|^2+|\llg_u(u)|^2\llg^{-1}(u)|Du|^2|u|^{2p}\\&\lesssim& \llg(u)|Du|^{2p}+\llg(u)|Du|^2|u|^{2p-2}+\llg(u)|u|^{2p}\lesssim \llg(u)|Du|^{2p}+\llg(u)|u|^{2p}.\earr$$ Therefore, by the assumptions \mref{piter1a} and \mref{IJp}, the last integral in \mref{Ipest} is bounded by a constant $C(C_0,R,p)$. We conclude that $I_{\cg_2p}(R)\le C(C_0,R,p)$. We repeat the argument with $g(u)=\llg(u)$ to see that $J_p(R)$ is also self improving. In this case $g(u)\llg^{-1}(u)\in L^\infty(Q)$ so that we can take $\cg_2$ to be any number in $(1,\cg_1)$. 

We let $\cg_*=\min\{\cg_1,\cg_2\}$ and complete the proof of the lemma. \eproof

\brem{murem1} It is also important to note that the estimate of \reflemm{dleANSpropt0}, based on those in  \reflemm{dleANSenergy}, \reflemm{dleANSprop}, is  {\em independent} of lower/upper bounds of the function $\llg_*$ in AR) but the integrals in M.0). The assumption AR) was used only in \reflemm{Tmaplem} to define the map $T_\sg$ and \reflemm{Dulocbound} to show that fixed points of $T_\sg$ are strong solutions. \erem

We are ready to provide the proof of the main theorem of this section. 

{\bf Proof of \reftheo{gentheo1}:} 
It is now clear that the assumptions M.0) and M.1) of our theorem allow us to
apply \reflemm{dleANSpropt0} and obtain a priori uniform bound  for any continuous strong solution $u$ of \mref{gensysfam}. The uniform estimate \mref{keydupANSt0} shows that the assumption \mref{M*defz} of \reflemm{Tmaplem} holds true so that the map $T_\sg$ is well defined and compact on a ball $\bar{B}_M$ of $\mX$ for some $M$ depending on the bound $M_*$ provided by \reflemm{dleANSpropt0}. Combining with \reflemm{Dulocbound}, the fixed points of $T_\sg$ are strong solutions of the system \mref{gensysfam} so that $T_\sg$ does not have a fixed point on the boundary of $\bar{B}_M$. Thus, by the Leray-Schauder fixed point theorem, $T_1$ has a fixed point in $B_M$ which is a strong solution to \mref{gensys}, which is unique because $u,Du$ are bounded and \mref{gensys} is now regular parabolic. The proof is complete.
\eproof

\section{Proof of the theorem on the general SKT system}\eqnoset\label{genSKTsec}

We conclude this paper by giving the proof of \reftheo{SKTgenthm}, an application of our main \reftheo{gentheo1}. To this end, we need only check the conditions A),H) and M.1) because the condition F) is obvious and M.0) is already assumed.

For $C^2$ positive scalar functions $\llg_i$, $i=1,\ldots,m$, and $\llg$ on $\RR^m$ we recall the notations in L): $\mL=\mbox{diag}[\llg_1(u),\ldots,\llg_m(u)],\,\mL_u=D_u[\llg_i(u)]_{i=1}^m$ and its assumptions
\beqno{llgcondz} \llg_i(u)\lesssim\llg(u),\; |u||\llg_u(u)|\lesssim \llg(u),\;|(\llg_i(u))_{uu}|\lesssim |\llg_{uu}(u)|.\eeq
\beqno{Aelliptic} \myprod{(\mL+\mbox{diag}[u_1,\ldots,u_m]\mL_u)\zeta,\zeta}\ge \llg(u)|\zeta|^2,\eeq
\beqno{mLhyp} |\mL_u|\lesssim |\llg_u(u)|,\; |\mL_u^{-1}|\lesssim |\llg_u(u)|^{-1}.\eeq

Recall that $P_i(u)=u_i\llg_i(u)$ so that
 $\partial_{u_j}P_i(u)=\dg_{ij}\llg_i(u)+u_i\partial_{u_j}\llg_i(u)$,  where $\dg_{ij}$ is the Kronecker delta. 
Writing $P(u)=[P_i(u)]_{i=1}^m$ and $\mU=\mbox{diag}[u_1,\ldots,u_m]$, we then have
$$A(u):=P_u(u)=[\partial_{u_j}P_i(u)]=\mL+\mU\mL_u.$$

We define $$K(u)=[K_i(u)]_{i=1}^m,\; K_i(u)=\log(\llg_i(u)).$$

We first have the following lemma.
\blemm{AHlem} The  matrix $A$ and the map $K$ satisfy the conditions A), H) respectively. \elemm
\bproof It is clear that \mref{Aelliptic} yields $\myprod{A(u)\zeta,\zeta}\ge \llg(u)|\zeta|^2$. The conditions \mref{llgcondz} and \mref{mLhyp} imply easily that $|A(u)|\lesssim\llg(u)$. Furthermore, simple calculation shows that they also give that $$\barr{lll}|A_u(u)|&\lesssim& |\mL_u|+|u||(\mL_u)u|\lesssim |\llg_u(u)|+|u|\max_i|(\llg_i(u))_{uu}|\\&\lesssim& |\llg_u(u)|+|u||\llg_u(u)|^2\llg^{-1}(u)\lesssim |\llg_u(u)|,\earr$$ because $|u||\llg_u(u)|\lesssim\llg(u)$. Thus, the condition A) is verified.

We turn to the map $K$. Because $\partial_{u_j}K_i(u)=\llg_i^{-1}(u)\partial_{u_j}\llg_i(u)$, we have $K_u(u)=\mL^{-1}\mL_u$. For $\mathbb{K}(u)=(K_u^{-1}(u))^T$ we have $|\mathbb{K}(u)|=|K_u^{-1}(u)|=|\mL_u^{-1}\mL|\lesssim \llg(u)|\llg_u(u)|^{-1}$ thanks to \mref{mLhyp}. On the other hand, $$|\mathbb{K}_u(u)|=|K_{uu}^{-1}(u)|=|(\mL_u^{-1}\mL)_u|\lesssim |\mL_u^{-1}||(\mL)_u|+\mL||(\mL_u)_u|||\mL_u^{-1}|^2.$$ 
Using the facts that $|(\mL_u)_u|\lesssim |\llg_{uu}(u)|$ and $\llg(u)|\llg_{uu}(u)|\lesssim|\llg_u(u)|^2$,  we see that $|\mathbb{K}_u(u)|$ is bounded by a constant.
Thus, H) is verified and the lemma is proved. \eproof

We now establish the first part of M.1) by showing that $K(u)$ is VMO. \blemm{SKT2d} Assume that there is a constant $C_0$ such that\beqno{LLubound} |\llg_u(u)|\llg^{-2}(u)\le C_0,\eeq \beqno{fSKT2d} \dspl{\int_0^{T_0}}\iidx{\Og}{|\llg(u)\bar{f}(u)|^2}\le C_0, \mbox{ where $\bar{f}(u):=[f_i(u)]_{i=1}^m$.}\eeq Then there is a constant $C(C_0)$ such that  \beqno{DKuest}\sup_{\tau\in(0,T_0)}\iidx{\Og\times\{\tau\}}{|D(K(u))|^2}\le C(C_0).\eeq 
\elemm

\bproof First of all, we test the $i$-th equation of the system with $(P_i(u))_t$ and sum the results to have, denoting $B(u,Du):=[B_i(u,Du)]_{i=1}^m$ \beqno{Ptest}\iidx{\Og}{\myprod{u_t,P_t(u)}}+\iidx{\Og}{\myprod{DP(u),D(P_t(u)}}=\iidx{\Og}{\myprod{B(u,Du)+\bar{f}(u),P_t(u)}}.\eeq 

Because $P_u(u)=A(u)$, $|A(u)|\lesssim \llg(u)$ and $|B(u,Du)|\lesssim\llg^\frac12(u)|Du|$, we have by Young's inequality
$$\myprod{B(u,Du)+\hat{f}(u),P_t(u)}\le \eg\llg(u)|u_t|^2 + C(\eg)\llg^2(u)|Du|^2+C(\eg)\llg(u)|\bar{f}(u)|^2.$$ As $\myprod{u_t,P_t(u)}\ge \llg(u)|u_t|^2$ and $\frac{d|DP(u)|^2}{dt}=\myprod{DP(u),D(P_t(u)}$ (because $u$ is a strong solution), we can choose $\eg$ small in the above and derive from \mref{Ptest} that \beqno{DPGronwall}\frac{d}{dt}\iidx{\Og}{|DP(u)|^2} \le C\iidx{\Og}{\llg^2(u)|Du|^2}+C\iidx{\Og}{\llg(u)|\bar{f}(u)|^2}.\eeq
By the ellipticity condition, $\llg(u)|Du|^2 \le \myprod{A(u),Du}\le \myprod{DP(u),Du}$ so that a simple use of Young's inequality implies $\llg(u)|Du|^2\le \frac12\llg^{-1}(u)|DP(u)|^2+\frac12\llg(u)|Du|^2$. This shows that $D(P(u))\sim \llg(u)|Du|$. We then have the following.  $$y'(t)\le Cy(t)+\ag(t), \mbox{ where } y(t)=\iidx{\Og\times\{t\}}{|DP(u)|^2},\; \ag(t)=\iidx{\Og\times\{t\}}{|\llg(u)\bar{f}(u)|^2}.$$ This is a simple Gronwall inequality for $y(t)$ and we have $$\sup_{\tau\in(0,T_0)}\iidx{\Og\times\{\tau\}}{|\llg(u)Du|^2}\sim\sup_{\tau\in(0,T_0)}\iidx{\Og\times\{\tau\}}{|DP(u)|^2}\le C\itQ{Q}{|\llg(u)\bar{f}(u)|^2}.$$

On the other hand, $|D(K(u))|=|K_u(u)Du|\lesssim |\mL^{-1}\mL_u|\llg^{-1}(u)|\llg(u)Du|$ so that  if  \mref{LLubound} holds then the above and \mref{fSKT2d} imply \mref{DKuest} and conclude the proof. \eproof

\brem{SKTgenBrem} The above lemma also shows that $K(u)$ is VMO($\Og$). We simply apply Poincar\'e's inequality for $n=2$ and use \mref{DKuest}. Also, the growth in $Du$ of $B_i(u,Du)$ is a bit different from f.1) in this paper but it was considered in \cite{dleANS} that $|\hat{f}(u,Du)|\lesssim \llg^{\frac12}(u)|Du|+f(u)$. \cite[Lemma 3.2]{dleANS} still provides the same energy estimate for $Du$ as \mref{keydupANSenergy}. Thus, the proof of our main theorem can continue.
\erem

{\bf Proof of \reftheo{SKTgenthm}:} By \reflemm{AHlem}, the assumption A), F) and H) of the main theorem are satisfied. As we already assumed M.0), the theorem will follow if M.1) is verified. The first part of M.1) requires that $K(u)$ has small BMO norm in small balls is given by \reflemm{SKT2d} and \refrem{SKTgenBrem}. We need only to check the second part by showing $\myPi_p(x,\tau):= \llg^{p+\frac12}(u)|\llg_u(u)|^{-p}$ is a weight. To this end, we will show that $\llg(u)$ and $|\llg_u(u)|$ are $A_1$ weights, $A_1=\cap_{\cg>1} A_\cg$. For $w_1=\log(\llg(u))$ and $w_2=\log(|\llg_u(u)|^{-1})$ we have $$|Dw_1|\le\frac{|\llg_u(u)|}{\llg(u)}|Du|\le \frac{|\llg_u(u)|}{\llg^2(u)}|\llg(u)Du|,$$ $$|Dw_2|\le\frac{|\llg_{uu}(u)|}{|\llg_u(u)|}|Du|\le\frac{|\llg_{uu}(u)|\llg(u)}{|\llg_u(u)|\llg^2(u)}|\llg(u)Du|\le \frac{|\llg_u(u)|}{\llg^2(u)}|\llg(u)Du|.$$ 

Since $|\llg_u(u)|/\llg^2(u)$ is bounded and $\llg(u)Du\in L^2(\Og)$, we see that $Dw_i\in L^2(\Og)$ and $w_i$'s have small BMO norm in small balls. It is wellknown that this implies $e^{cw_i}$, and therefore $\llg^c(u)$ and $|\llg_u(u)|^{-c}$, are $A_1$ weights for any $c>0$ (see \cite{Graf} or \cite[Lemma 5.1]{dleAMJ}). Hence, for each $\tau\in(0,T_0)$ and any power of $\myPi_p(x,\tau)$ is also an $A_1$ weight and the last condition in H) is then verified. The proof is complete. \eproof

\bibliographystyle{plain}

\end{document}